\documentclass[preprint,groupedaddress,amsmath,amssymb,showkeys]{revtex4}

\usepackage{dcolumn}
\usepackage{bm}

\usepackage{amsmath}
\usepackage{amssymb}
\usepackage{latexsym}
\usepackage[dvips]{graphicx}
\usepackage{afterpage}
\usepackage{float}
\newcommand{\m}{\mbox}
\newcommand{\bd}{\boldmath}
\newcommand{\bo}[1]{\mbox{\boldmath $#1$}}

\begin{document}


\title{Magnetically-induced buckling of a whirling conducting rod with applications to electrodynamic space tethers}


\author{J.~Valverde}
\affiliation{Department of Civil and Environmental Engineering, University
of California, Berkeley, USA}
\author{G.H.M.~van~der~Heijden}
\email{g.heijden@ucl.ac.uk}
\affiliation{Centre for Nonlinear Dynamics, University College London,\\
 Gower Street, London WC1E 6BT, UK}

\date{\today}

\begin{abstract}
We study the effect of a magnetic field on the behaviour of a
slender conducting elastic structure, motivated by stability problems
of electrodynamic space tethers. Both statical (buckling) and dynamical
(whirling) instability are considered and we also compute post-buckling
configurations. The equations used are the geometrically exact Kirchhoff
equations. Magnetic buckling of a welded rod is found to be described by
a surprisingly degenerate bifurcation, which is unfolded when both
transverse anisotropy of the rod and angular velocity are considered.
By solving the linearised equations about the (quasi-) stationary
solutions we find various secondary instabilities. Our results are relevant
for current designs of electrodynamic space tethers and potentially for
future applications in nano- and molecular wires.
\end{abstract}


\maketitle


\section{Introduction}

\begin{figure}
\begin{center}
\includegraphics[width=0.3\linewidth]{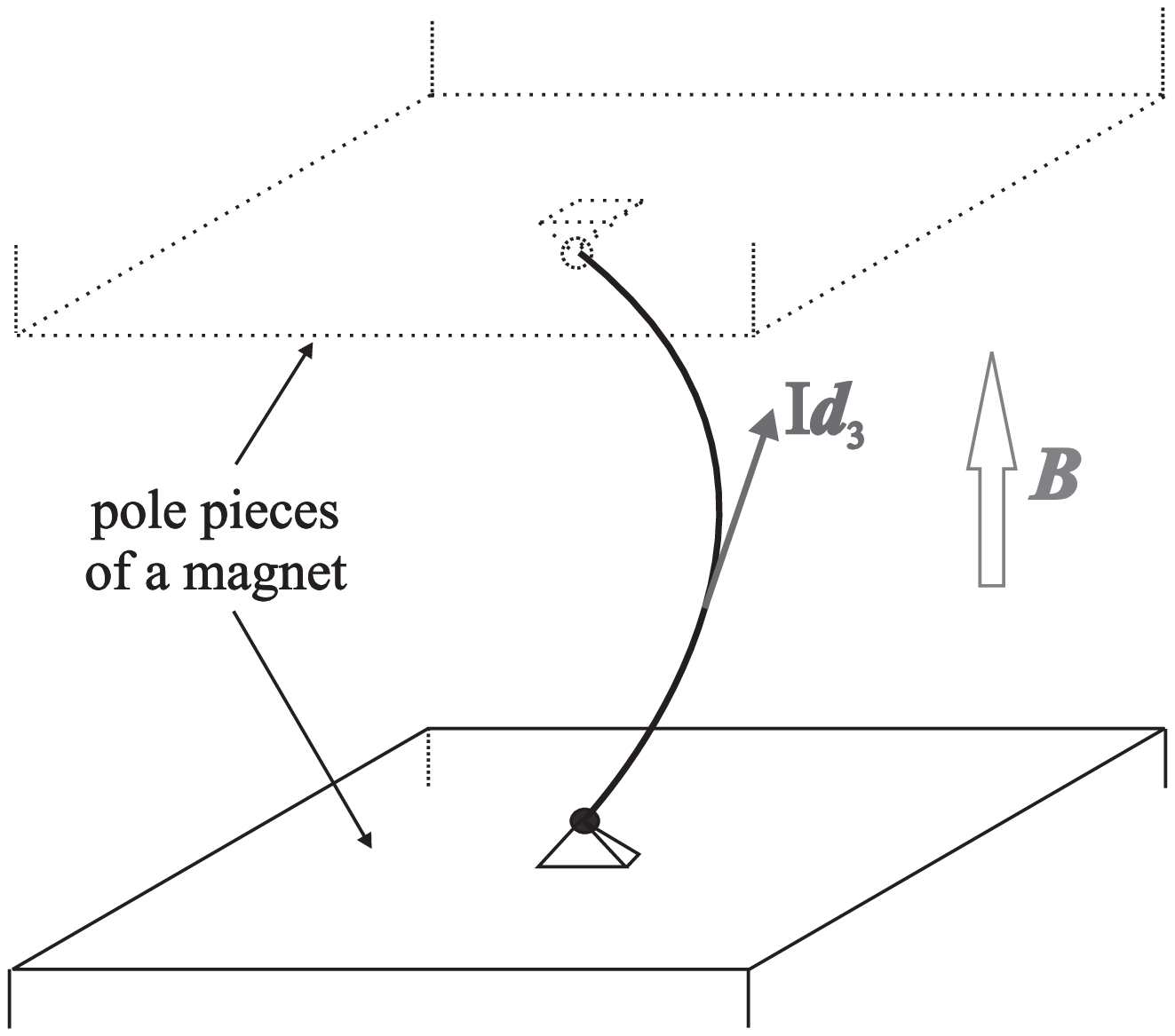}
\end{center}
\caption{Experimental setup for a conducting wire.}
\label{fig:hel:exp-setup}
\end{figure}

A straight current-carrying wire held in tension between pole faces
of a magnet is well known to buckle into a (roughly) helical
configuration at a critical current (see
Figure~\ref{fig:hel:exp-setup}). A photograph of this phenomenon is
shown in Section 10.4.3 of \cite{woodson}, where
a linear stability analysis is carried out for a simple string
model. (A string is here meant to be a perfectly flexible elastic
wire.) The problem was studied by Wolfe \cite{wolfe1} by means of a
rigorous bifurcation analysis for a (nonlinearly elastic) string
suspended between fixed supports and placed in a uniform magnetic
field directed parallel to the undeformed wire. He found that an
infinite number of solution branches bifurcate from the trivial
straight solution, much like in the Euler elastica under compressive
load. Minimisation of the potential energy indicated that the first
branch of bifurcating solutions is stable while all other branches are
unstable.

In a subsequent paper Wolfe \cite{wolfe2} extends the analysis to a
uniformly rotating (whirling) string and shows again the existence
of bifurcating branches of whirling non-trivial solutions. This result
was further extended by Healey \cite{healey} using equivariant
bifurcation theory in order to deal with the symmetries of the problem,
which caused the bifurcations to be degenerate.

Wolfe also considered a conducting {\it rod} in a uniform magnetic
field \cite{wolfe3}. In addition to extension a rod can undergo
flexure, torsion and shear, and for the case of welded boundary
conditions it was found that in certain cases bifurcation occurs,
with the usual infinity of non-trivial equilibrium states. All the
works cited above were content with showing the existence of
bifurcating solutions and did not study their post-buckling
behaviour. The Hamiltonian structure of the equations for a rod in
a magnetic field was investigated in \cite{SH} where it was shown
that in the case of an isotropic, inextensible and unshearable rod
the equations are completely integrable.

The study of strings and rods in a magnetic field is of great
interest to space tethers. Although space tethers in the last 20
years or so have become a well established concept in astrodynamics
\cite{beletsky}, new designs continue to be proposed that hold great
potential for future space applications. A space tether is a long
cable used to connect spacecraft to each other or to other orbiting
bodies such as space stations, boosters, payloads, etc.~in order to
transfer energy and momentum thus providing space propulsion without
consuming propellant. These tethers have been studied as elastic
strings (e.g., \cite{krupa_et_al}) and as dumb-bell systems (e.g.,
\cite{ziegler_cartmell}). An important class of space tethers is
formed by the so-called electrodynamic tethers (ETs). These employ
the earth's magnetic field and ionospheric plasma to generate a
current, according to Faraday's Law, and hence thrust or drag forces
without expending chemical fuel. An example is the Short
Electrodynamic Tether (SET) prototype of the European Space Agency
\cite{valverde,valverde-cosserat}. This tether system, which
comprises a central module from which two tethers each about a
hundred metres long extend, is designed to operate at high
inclination and in low orbit. Due to the shape of the earth's
magnetic field this means that the desirable orientation of the
tether is the horizontal one (i.e., with the axis of smallest moment
of inertia normal to the orbit plane). The gravity gradient
ordinarily causes a tether to drift to the stable radial position.
Therefore, in order to keep the system in the horizontal position,
an axial spin velocity is applied for gyroscopic stability
\cite{sanmartin,hughes_stability,valverde}. This requires
significant torsional and bending rigidity of the tether which
therefore has to be modelled as a rod, not a string. The applied
spin causes large deformations that present stability issues similar
to those in an unbalanced rotor system, which have been studied in
previous work \cite{valverde-cosserat,valverde2}.

In this paper we apply large-deformation rod theory to study the
effect of a magnetic field on the behaviour of a slender conducting
elastic structure possibly subject to end forces. Both statical
(buckling) and dynamical (whirling) instability are considered and
we also compute post-buckling configurations. The work extends the
stability analysis of the SET in \cite{valverde-cosserat} by
including the effect of the magnetic field. This effect was
considered small for the SET but that need not be true for other
tether designs. For instance, very long and flexible tethers subject
to boundary conditions that are not too restraining (e.g., no big
end masses) might well operate in the region of the first magnetic
buckling instability.

We consider welded boundary conditions. These are appropriate for tethers
with sufficiently large attached end masses (see
Figure~\ref{fig:et:mag-force}). Unlike in string buckling the rod does not
require a tensile force in the trivial state, but we allow for such an
applied force as well. The pertinent dimensionless parameter that governs
buckling measures the product of current and magnetic field against the
bending force. We find that magnetic buckling of the welded rod is described
by a remarkably degenerate pitchfork bifurcation. Wolfe considered welded
boundary conditions in his statical study in \cite{wolfe3} and encountered
degeneracies (even-dimensional eigenspaces) because of rotational symmetry
of the problem, but we show that further complications occur, involving
branches connecting the bifurcating branches.

We also study steady whirling solutions for which we introduce a
rotating coordinate system. This extends Wolfe's analysis of
whirling strings to whirling rods. It is found that applied spin
resolves the degeneracies of the pitchfork bifurcations, provided that
the rod is transversely anisotropic. We perform a stability analysis by
computing eigenvalues of the linearised boudary-value problem about a
(quasi-) stationary solution. For this we use a continuation (or homotopy)
approach that takes advantage of the fact that exact expressions for the
(imaginary) eigenvalues can be obtained in an appropriate limit (no spin,
no magnetic field). The eigenvalues in this limit are then traced as
system parameters are varied.

The paper is organised as follows. First, in Section 2, we give more
details about the tether application and discuss the effect of the
earth's magnetic field on an electrodynamic tether. Then we present
our rod mechanics formulation in Section 3, in which the magnetic
field enters the force balance equation through the Lorentz body
force. For the study of whirling solutions the equilibrium equations
are transformed to a coordinate system rotating at constant angular
velocity (Section 4). The linearisation is presented in Section 5 and our
continuation approach to stability analysis is discussed in Section 6.
Results are presented in the form of bifurcation diagrams in Section 7,
both for the statical and dynamical case. We also pick up some secondary
bifurcations with associated loss of stability. Finally, in Section 8 we
comment on the implications of our results for electrodynamic tethers and
draw some conclusions.

%
%
%
%
%

\section{Electrodynamic tethers and the Earth's magnetic field}

Electrodynamic tethers are electrical conductors that interact with
the geomagnetic field in such a way that an electromotive force
(e.m.f.) is generated along the tether due to Faraday's Law
\cite{lorenzini}. The electrical circuit is closed by means of two
contactors attached to the ends of the tether which interact with
the surrounding ionospheric plasma and allow a current to flow.


\begin{figure}
\begin{center}
\includegraphics[width=0.2\linewidth]{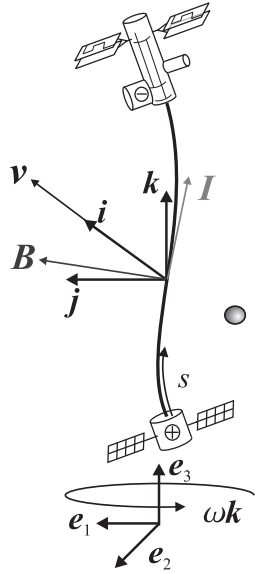}
\end{center}
\caption{An electrodynamic tether in the Earth's magnetic field.}
\label{fig:et:mag-force}
\end{figure}

Figure~\ref{fig:et:mag-force} shows an electrodynamic tether
connecting two satellites. Let $\{\bo i,\bo j,\bo k\}$ be the
orbital frame, which may be assumed to be inertial. This is due to
the fact that external forces on the tether system, that is,
electrodynamic forces, are small enough to consider a circular
Keplerian orbit which allows to decouple translational orbital
dynamics from structural deformation dynamics~\cite{lorenzini}.
Besides, since the angular velocity of the tether around the earth
is much smaller than the typical spin velocity about its axis,
rotational orbital dynamics can also be decoupled from structural
dynamics of the tether subjected to spin rotation, which is the main
issue of the present paper~\cite{valverde-cosserat}. In such a case,
the tether is assumed to travel with a velocity $\bo v$ in the
direction of $\bo i$. The e.m.f. between the ends induced by this
motion is given by
\begin{equation}
E = \int_l (\bo v \times \bo B_0) \cdot d \bo l, \label{eq:faraday}
\end{equation}
where $\bo B_0$ is the magnetic field and $d \bo l$ is a
differential along the length of the tether. Because the tether is
part of a closed circuit a current $\bo I$ will flow in the
direction of increasing $E$ and the system functions as a generator.
This current in turn gives rise to a Lorentz force ${\bo F}_L$
through
\begin{equation}
d {\bo F}_L = d \bo I \times \bo B_0. \label{eq:lorentz}
\end{equation}
This force can be used to drag the system without expending chemical
fuel \cite{lorenzini,beletsky}. Alternatively, if a current is
forced against the e.m.f. the system becomes a motor boosting itself
to a higher orbit.

If we denote the position co-ordinates of the tether's central axis
relative to $\{\bo i,\bo j,\bo k\}$ by $(X,Y,Z)$, then the current
vector, which is directed along the tangent of the tether, can be
expressed as
\begin{equation}
d \bo I = I d \bo l = I\left(\frac{\partial X}{\partial s},
\frac{\partial Y}{\partial s},\frac{\partial Z}{\partial s}\right)^T
ds. \label{eq:current}
\end{equation}
The maximum force is generated when $\bo I$ and $\bo B_0$ are
perpendicular. In the ET operation conditions both vectors will in
general not be perpendicular because the tether is not perfectly
straight and the magnetic field lines will not be perpendicular to
the tether over its entire length. In order to represent this
imperfection let us therefore assume that the magnetic field has an
extra component in the $\bo k$ direction,
\begin{equation}
\bo B_0= (0,B_1,B_2)^T, \label{eq:magnetic}
\end{equation}

\noindent where the desired component of the field $B_1 \gg B_2$.
Introducing (\ref{eq:current}) and (\ref{eq:magnetic}) into
(\ref{eq:lorentz}), the differential of the Lorentz force is found
to be
%
\begin{equation}
\textrm{d} \bo F_L =  I\left(\frac{\partial Y}{\partial s}
B_2-\frac{\partial Z}{\partial s} B_1,-\frac{\partial X} {\partial
s} B_2,\frac{\partial X}{\partial s} B_1\right)^T \textrm{d}s.
\label{eq:lorentz2}
\end{equation}

\noindent The term $- \frac{\partial Z}{\partial s} B_1 I ds$ in the
$\bo i$ component opposes the motion and drags the ET system, as
intended. However, the crossed $B_2$ terms in the $\bo i$ and $\bo
j$ components (due to the imperfection $B_2$) will tend to coil the
tether. This undesirable effect has been reported in some tether
flights~\cite{lorenzini,beletsky}.

In this paper we study the interaction of elastic and
electromagnetic forces in a conducting rod and the tendency to
generate three-dimensional coiled configurations (for instance as a
result of buckling at critical loads).


\section{The rod mechanics model}

We describe the elastic behaviour of a conducting cable by the
Kirchhoff equations for the dynamics of thin rods. The rod is
assumed to be uniform, inextensible, unshearable and intrinsically
straight and prismatic. The inextensible and unshearable assumptions
are appropriate for thin rods with low external forces
(electrodynamic), as the tether. For the background of the
Kirchhoff equations the reader is referred to
\cite{antman,coleman_et_al}. These equations were also used in
\cite{valverde-cosserat} to analyse the dynamics of the SET. To
these equations will be added the Lorentz force to account for the
electromagnetic interaction.

\begin{figure}
\begin{center}
\includegraphics[width=0.5\linewidth]{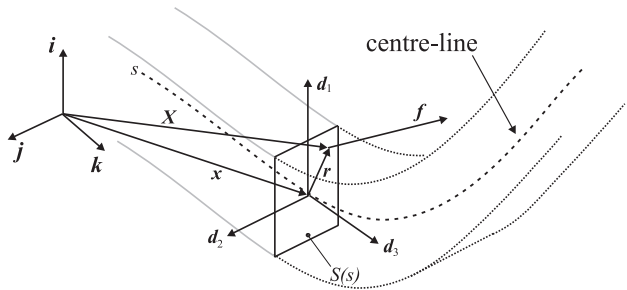}
\end{center}
\caption{Cosserat model of a rod.} \label{fig:cosserat}
\end{figure}

Let $\bo x$ denote the position of the rod's centreline and let
$\{\bo d_1,\bo d_2,\bo d_3\}$ be a right-handed orthonormal frame of
directors (the Cosserat triad) defined at each point along the
centreline. Since the centreline is assumed to be inextensible we
can take $\bo d_3$ in the direction of the local tangent:
\begin{equation}
\bo x'(s,t)=\bo d_3(s,t), \label{eqtangent1}
\end{equation}

\noindent where the prime denotes differentiation with respect to
arclength $s$ measured along the centreline, and $t$ is time. The
directors $\bo d_1$ and $\bo d_2$ will be taken to point along the
principal bending axes of the cross-section (see
Figure~\ref{fig:cosserat}). The unstressed rod is taken to lie along
the basis vector $\bo k$ of a fixed inertial frame $\{\bo i,\bo
j,\bo k\}$.

Looking at Figure~\ref{fig:cosserat} we note that the position
vector of an arbitrary point of the rod can be expressed as
\begin{equation}
\begin{split}
\bo X (s,\xi_1,\xi_2,t)=&\bo x (s,t) +\xi_1\bo d _1(s,t) +\xi_2\bo d
_2(s,t)\\
=&\bo x (s,t) +\bo r (s,\xi_1,\xi_2,t), \label{eq:kinem}
\end{split}
\end{equation}

\noindent where $(\xi_1,\xi_2)$ are the components of $\bo r$ in the
cross-section relative to $\{\m{\bd$d$}_1(s),\m{\bd$d$}_2(s)\}$. The
rod is thus viewed as a set of infinitesimal slices centred at all
$s$. A one-dimensional description will be obtained by averaging of
forces and moments over each cross-section. The internal traction,
which is the projection of the stress tensor onto the
cross-sectional plane, is given by a force which we denote by $\bo f
=\bo f (s,\xi_1,\xi_2,t)$ (see Figure~\ref{fig:cosserat}). The
resultant elastic force exerted in a section $S(s)$ is given by
\begin{equation}
\bo F (s,t)=\int_{S(s)} \bo f (s,\xi_1,\xi_2,t)\,\textrm{d}S,
\end{equation}
\noindent where $\textrm{d}S$ is an infinitesimal area element. This
force can be expressed in the director basis as $\bo F =
\sum_{i=1}^{3} F_{i}\,\m{\bd$d$}_i$. The resultant moment in the
section $S(s)$ is given by
\begin{equation}
\bo M (s,t)=\int_{S(s)} \bo r (s,\xi_1,\xi_2,t) \times  \bo f
(s,\xi_1,\xi_2,t)\,\textrm{d}S, \label{moments}
\end{equation}
\noindent and will be expressed as $\bo M = \sum_{i=1}^{3}
M_{i}\m{\bd$d$}_i$.

The rod is assumed to carry an electric current for which we can
write
\begin{equation}
\bo I = I \bo x' = I \bo d_3.
\end{equation}
Here we have assumed the current to have the same direction as the
rod, which is consistent with a one-dimensional rod theory. It
amounts to the assumption that the cross-section of the conducting
wire is small enough to make currents within the cross-section (eddy
currents) induced by the motion negligible. The current $\bo I$
interacts with the magnetic field $\bo B_0$ to generate a (Lorentz)
body force given by
\begin{equation}
\bo F_L= I\bo d_3 \times \bo B_0.
\end{equation}
Following \cite{wolfe1} we assume the magnetic field to be uniform
and directed along the unstressed rod, i.e.,
\begin{equation}
\bo B_0= B_0 \bo k,
\end{equation}
modelling the undesired component of the magnetic field, $B_2$ in
equation (\ref{eq:magnetic}).

The balancing of forces and moments across an infinitesimal rod
element~\cite{valverde-cosserat} then yields the following set of
partial differential equations:
\begin{equation}
\bo F' + IB_0 \bo d_3 \times \bo k = \rho A  \ddot{\bo x} ,
\label{eq:lmomentum2}
\end{equation}
\begin{equation}
\bo M'  + \bo d_3  \times \bo F= \rho ( I_2 \bo d_1 \times \ddot{\bo
d_1}+ I_1 \bo d_2 \times \ddot{\bo d_2}), \label{eq:amomentum2}
\end{equation}

\noindent where $\rho$ is the (volumetric) mass density, $A$ the
cross-sectional area, $I_1$ and $I_2$ the second moment of area of
the cross-section about $\bo d_1$ and $\bo d_2$ respectively, and
$\dot{(~)}$ denotes differentiation with respect to time.


For a closed system of equations these balance equations need to be
supplemented by constitutive relations that characterise the
material behaviour of the rod. We assume the rod to be made of
homogeneous isotropic linear viscoelastic material so that
stress-strain relations, based on a model by Valverde et
al.~\cite{valverde-cosserat}, are
\begin{equation}
\begin{split}
M_1 = EI_1 (\kappa_1+\gamma_v \dot \kappa_1),\\
M_2 = EI_2 (\kappa_2+\gamma_v \dot \kappa_2), \\
M_3 = GJ (\kappa_3+\gamma_v \dot \kappa_3), \\
\end{split}
\label{eq:constitutive1}
\end{equation}
where $\kappa_1$ and $\kappa_2$ are the curvatures about $\bo d_1$
and $\bo d_2$, respectively, while $\kappa_3$ is the twist about
$\bo d_3$. The constant $\gamma_v$ is the viscoelastic coefficient
of the material, $E$ is Young's modulus, $G$ is the shear modulus
and $J$ is the second moment of area of the cross-section about $\bo
d_3$. We shall assume that the section is symmetric with respect to
the principal axes, in which case $J=I_1+I_2$.

The $\kappa_i$ are the components of the curvature vector
\begin{equation}
\m{\bd$\kappa$}=\sum_{i=1}^{3} \kappa_{i}\,\m{\bd$d$}_i,
\label{eq:twist1}
\end{equation}
which governs the evolution in space of the frame of directors as
one moves along the centreline:
\begin{equation}
\m{\bd$d$}'_i=\m{\bd$\kappa$}\times \m{\bd$d$}_i \quad\quad
(i=1,2,3). \label{eq:twist}
\end{equation}

The constitutive relations (\ref{eq:constitutive1}) can be used to
replace the $\kappa_i$ in (\ref{eq:twist}) by moments, after which
the equations (\ref{eqtangent1}), (\ref{eq:lmomentum2}),
(\ref{eq:amomentum2}) and (\ref{eq:twist}) form a system of 18
differential equations for the 18
unknowns $(\bo x,\bo F,\bo M,\bo d_1,\bo d_2,\bo d_3)$.\\

\noindent {\bf Remark:} In general when a conducting wire moves in a
magnetic field an additional electromagnetic induction effect occurs
which opposes the motion. The electromotive force as a result of
this effect is proportional to the rate of change of the enclosed
magnetic flux \cite{jackson}. However, in the case of a steadily
whirling wire in a uniform magnetic field the enclosed magnetic flux
does not change and the effect is zero.

\subsection{Equations of motion in the moving frame}
Since we are interested in steadily rotating solutions we transform
the equilibrium equations (\ref{eq:lmomentum2}) and
(\ref{eq:amomentum2}) to a co-ordinate frame $\{\bo e_1,\bo e_2,\bo
e_3\}$ that rotates with constant angular velocity $\bo
\omega=\omega \bo k$ about the $\bo k$ axis (and the axis of the rod
in its trivial unstressed state). Noting that the derivative with
respect to time of an arbitrary vector $\bo V(s,t)$ is given by
\begin{equation}
\left . \frac {d \bo V (s,t)}{d t}\right|_i=\left . \frac {d \bo V
(s,t)}{d t}\right|_m + \bo \omega \times \bo V (s,t),
\label{eq:derfixed}
\end{equation}

\noindent where $\left . \frac {d}{d t}\right|_i$ indicates the derivative
with respect to time in the inertial frame and $\left . \frac {d}{d
t}\right|_m$ stands for the derivative with respect to time in the
moving frame, the equations (\ref{eq:lmomentum2}) and
(\ref{eq:amomentum2}) expressed relative to $\{\bo e_1,\bo e_2,\bo
e_3\}$ become
\begin{equation}
\bo F' + IB_0\bo d_3 \times \bo e_3 = \rho A (\ddot{\bo x} +2 \bo
\omega \times \dot{\bo x} + \bo \omega \times (\bo \omega \times \bo
x)), \label{eq:lmomentum4}
\end{equation}
\begin{equation}
\begin{split}
\bo M'  + \bo d_3  \times \bo F &= \rho I_2 (\bo d_1 \times
\ddot{\bo d_1}+ 2\bo d_1 \times (\bo \omega \times \dot{\bo
d_1})+(\bo \omega \cdot
\bo d_1)(\bo d_1 \times \bo \omega)) \\
&+ \rho I_1 (\bo d_2 \times \ddot{\bo d_2} + 2\bo d_2 \times (\bo
\omega \times \dot{\bo d_2})+(\bo \omega \cdot \bo d_2)(\bo d_2
\times \bo \omega)). \label{eq:amomentum4}
\end{split}
\end{equation}

Steadily rotating (whirling) solutions satisfy the equations
(\ref{eq:lmomentum4}) and (\ref{eq:amomentum4}) with the dotted
variables set to zero:
\begin{equation}
\bo F' + IB_0\bo d_3 \times \bo e_3 = \rho A \bo \omega \times (\bo
\omega \times \bo x), \label{eq:lmomentum4_whirl}
\end{equation}
\begin{equation}
\bo M'  + \bo d_3  \times \bo F = \rho I_2 (\bo \omega \cdot\bo d_1)
(\bo d_1 \times \bo \omega) + \rho I_1 (\bo \omega \cdot \bo
d_2)(\bo d_2\times \bo \omega). \label{eq:amomentum4_whirl}
\end{equation}
The other equations (\ref{eqtangent1}) and (\ref{eq:twist}) do not
change their form, but all vectors are now to be considered as
expressed relative to the rotating frame $\{\bo e_1,\bo e_2,\bo
e_3\}$. Statical solutions are simply obtained by setting $\omega$
equal to zero.

For a well-posed problem the final 18 ODEs require 18 boundary
conditions to be specified.

\subsection{The boundary conditions}

We follow \cite{wolfe3} and consider welded boundary conditions.
These conditions also describe an electrodynamic tether that is
welded to the end contactors or modules if these bodies are
sufficiently massive (see Figure~\ref{fig:et:mag-force}). We assume
the rod to be fixed at $s=L$ and to be able to slide along $\bo
e_3=\bo k$ at $s=0$ where a controlled force $T$ is applied
(positive for tension). $L$ is the length of the rod. Writing $\bo
x=x\bo e_1 + y\bo e_2 + z\bo e_3$ a consistent set of boundary
conditions is thus given by
{\setlength\arraycolsep{2pt}
\begin{eqnarray}
x (0,t)= 0, \label{eq:et:BC23}\\
y (0,t)= 0, \label{eq:et:BC24}\\
\bo F (0,t) \cdot \bo e_3 = T, \label{eq:et:BC25}\\
\bo d_3 (0,t) \cdot \bo e_1 = 0, \label{eq:et:BC26}\\
\bo d_3 (0,t)\cdot \bo e_2 = 0, \label{eq:et:BC27}\\
\bo d_1 (0,t) \cdot \bo e_2 = 0, \label{eq:et:BC28}
\end{eqnarray}}
at $s=0$ and
{\setlength\arraycolsep{2pt}
\begin{eqnarray}
x (L,t)= 0, \label{eq:et:BC29}\\
y (L,t)= 0, \label{eq:et:BC30}\\
z (L,t)= L, \label{eq:et:BC31}\\
\bo d_3 (L,t) \cdot \bo e_1 = 0, \label{eq:et:BC32}\\
\bo d_3 (L,t) \cdot \bo e_2 = 0, \label{eq:et:BC33}\\
\bo d_1 (L,t) \cdot \bo e_2 = 0, \label{eq:et:BC34}
\end{eqnarray}}

\noindent at $s=L$. To these conditions we have to add conditions
that ensure the orthonormality of the director basis, for which we
can take
\begin{equation}
\begin{split}
\bo d_1 (0,t) \cdot \bo d_1 (0,t) = 1, \\
\bo d_2 (0,t) \cdot \bo d_2 (0,t) = 1, \\
\bo d_3 (0,t) \cdot \bo d_3 (0,t) = 1, \\
\bo d_1 (0,t) \cdot \bo d_2 (0,t) = 0, \\
\bo d_1 (0,t) \cdot \bo d_3 (0,t) = 0, \\
\bo d_2 (0,t) \cdot \bo d_3 (0,t) = 0, \\
\end{split}
\label{eq:et:orthonorm}
\end{equation}
for a total of 18 boundary conditions, as required. Note that these
conditions imply that at $s=0$ and $s=L$ the director frame $\{\bo
d_1,\bo d_2,\bo d_3\}$ is aligned with $\{\bo e_1,\bo e_2,\bo
e_3\}$.

\section{Nondimensionalisation}

We make the system of equations dimensionless by scaling the
variables in the following way
\begin{equation}
\begin{split}
&\omega_c = \sqrt{f \frac{EI_1}{\rho A L^4}},\quad \bar t=t\omega_c,
\quad \bar s=\frac{s}{L}\in[0,1], \quad \bar \omega=\frac{\omega}{\omega_c},
\quad \bar{\bo x}=\frac{\bo x}{L},
\\
&\bar {\bo F}=\bo F \frac{L^2}{f EI_1},\quad \bar{T}=T \frac{L^2}{f E I_1},
\quad \bar{\bo M}=\bo M\frac{L}{f E I_1}, \quad \bar{\bo \kappa}=
\bo\kappa L.
\end{split}
\label{eq:nondim}
\end{equation}

\noindent Here $\omega_c$ is a reference bending natural frequency
which, through tuning of the numerical constant $f$, can be adapted
to the particular boundary conditions at hand and the natural mode
considered. In the analysis in Section \ref{sect:stab_anal} we shall take
$f=1$, but since in the welded case we consider realistic data for
electrodynamic tethers we shall take $f=500.5639$ in that case so as to
get convenient numbers when presenting our numerical results. This value
for $f$ corresponds to the first bending natural frequency of a
welded-welded beam about $\bo d_1$.

With this nondimensionalisation the equations become (dropping the
overbars for simplicity and letting a prime denote $\frac{d}{d \bar s}$):
\begin{eqnarray}
&\bo F' + B \bo d_3 \times \bo k =\ddot{\bo x} +2 \bo \omega \times
\dot{\bo x} + \bo \omega \times (\bo \omega
\times \bo x), \label{eq:lmomentum}\\
\bo M' + \bo d_3  \times \bo F &= P\left[R (\bo d_1 \times \ddot{\bo
d_1}+ 2\bo d_1 \times \bo \omega \times \dot{\bo d_1}+(\bo \omega
\cdot \bo d_1)(\bo d_1 \times \bo \omega)) \right. \nonumber \\
& \left. + (\bo d_2 \times \ddot{\bo d_2} + 2\bo d_2 \times \bo \omega \times
\dot{\bo d_2}+(\bo \omega \cdot \bo d_2)(\bo d_2 \times \bo
\omega)) \right], \label{eq:amomentum}\\
&\bo x'=\bo d_3, \label{eq:tangent}\\
&\bo d_i'=\bo \kappa \times \bo d_i, \label{eq:twist2}
\end{eqnarray}
\noindent and the constitutive relations can be written as
\begin{equation}
\bo M=\frac{1}{f}\left[(\kappa_1+\gamma \dot{\kappa_1})\bo d_1+
R(\kappa_2+\gamma \dot{\kappa_2})\bo d_2+ \frac{\Gamma
(1+R)}{2}(\kappa_3+\gamma \dot{\kappa_3})\bo d_3\right],
\label{eq:constitutive}
\end{equation}
\noindent where the dimensionless parameters are
\begin{equation}
P=\frac{I_1}{A L^2},\quad R=\frac{I_2}{I_1},\quad B=\frac{B_0 I
L^3}{f E I_1},\quad \Gamma = \frac{2G}{E},\quad
\gamma=\gamma_v\omega_c,
\end{equation}
and ($\frac{1}{\Gamma}-1$) is equal to Poisson's ratio. For the
boundary conditions we can still use (\ref{eq:et:BC23}) to
(\ref{eq:et:orthonorm}) if we assume that they now refer to
dimensionless variables and that the right-hand conditions are
imposed at $\bar{s}=1$.

\section{Perturbation scheme}

We consider whirling solutions (relative equilibria) that are stationary
in the moving frame $\{\bo e_1,\bo e_2,\bo e_3\}$. Such solutions are
found by solving the set of equations
(\ref{eq:lmomentum})--(\ref{eq:constitutive}) with the dotted variables set
to zero (thus obtaining an ODE). To study their stability we linearise the
full PDE (\ref{eq:lmomentum})--(\ref{eq:constitutive}) about these whirling
solutions. We follow the approach in \cite{valverde-cosserat}, which is
similar to approaches in \cite{goriely,fraser}. The stability of statical
(non-whirling) solutions can be investigated by simply setting the angular
velocity $\omega$ to zero.

We start our perturbation analysis by writing
\begin{equation}
\bo d_i(s,t)=\bo d_i^{0}(s)+\delta \bo
d_i^{t}(s,t)+O(\delta^2),\quad i=1,2,3, \label{per:dir}
\end{equation}
where $\bo d_i^{0}(s)$ are the basis vectors of a quasi-stationary solution,
$\bo d_i^{t}(s,t)$ are the basis vectors of a time-dependent perturbation
and $\delta$ is a small bookkeeping parameter introduced to separate scales.
Note that, in order to preserve orthonormality to $O(\delta)$ ($\bo d_i
\cdot \bo d_j = \delta_{ij}+O(\delta^2)$), we must have
\begin{equation}
\bo d_i^{t}(s,t)=\sum_{j=1}^3 A_{ij}(s,t) \bo d_j^0(s),\quad
i=1,2,3, \label{per:dirt}
\end{equation}
where the matrix $A_{ij}$ is skew-symmetric and can be written as
\begin{equation}
\bo A= \left( \begin{array}{ccc}
0 & \alpha_3 & -\alpha_2 \\
-\alpha_3 & 0 & \alpha_1 \\
\alpha_2 & -\alpha_1  & 0
\end{array} \right).
\label{per:dirt2}
\end{equation}
Thus, the nine components of the director basis perturbation are
described by only three independent parameters, and if we introduce
\begin{equation}
\bo \alpha = (\alpha_1,\alpha_2,\alpha_3)^T, \label{alpha}
\end{equation}
(with respect to the unperturbed director basis) then the perturbed director
basis can be expressed as
\begin{equation}
\bo d_i(s,t)=\bo d_i^{0}(s)+\delta \bo \alpha(s,t) \times \bo d_i^{0}(s)
+O(\delta^2), \quad i=1,2,3. \label{per:dirfinal}
\end{equation}
%

Using (\ref{per:dirfinal}), the perturbation of an arbitrary vector
$\bo V=\sum_{i=1}^3 V_i \bo d_i$ can be written on the basis
$\{\bo d_1^0,\bo d_2^0,\bo d_3^0\}$ as
\begin{equation}
\bo V=\bo V^0+\delta \bo V^t  +O(\delta^2)
=\sum_{i=1}^3 [V_i^0+\delta( V_i^t +(\bo \alpha \times \bo
V^0)_i)] \bo d_i^0    +O(\delta^2), \label{per:varfinal}
\end{equation}
where $()_i$ denotes the component along $\bo d_i^0$ and time and space
dependence of the variables have been suppressed for the sake of simplicity.

Applying this perturbation scheme to the PDEs
(\ref{eq:lmomentum})--(\ref{eq:constitutive}) and the boundary conditions,
we arrive at an $O(1)$ nonlinear ODE for the quasi-stationary solutions and
an $O(\delta)$ linear PDE governing their stability.

\subsection{The $O(1)$ equations -- quasi-stationary whirl}
\label{stationaryBVP}

The $O(1)$ equations are time-independent. Recalling that $\bo
\omega= \omega \bo e_3$, we find the $O(1)$ terms of the linear
momentum equation (\ref{eq:lmomentum}), projected on the director
basis $\{\bo d_1^0,\bo d_2^0,\bo d_3^0\}$, to give

{\setlength\arraycolsep{2pt}
\begin{eqnarray}
(F_1^0)'-F_2^0\kappa_3^0+F_3^0\kappa_2^0 + B(d_{32}^0d_{11}^0-d_{31}^0
d_{12}^0)&=&-\omega^2(x^0d_{11}^0+y^0d_{12}^0),\label{stat:forces1}\\
(F_2^0)'-F_3^0\kappa_1^0+F_1^0\kappa_3^0 + B(d_{32}^0d_{21}^0-d_{31}^0
d_{22}^0)&=&-\omega^2(x^0d_{21}^0+y^0d_{22}^0),\label{stat:forces2}\\
(F_3^0)'-F_1^0\kappa_2^0+F_2^0\kappa_1^0&=&-\omega^2(x^0d_{31}^0+y^0d_{32}^0),
\label{stat:forces3}
\end{eqnarray}}

\noindent where subscripts are used to indicate components relative
to the basis vectors $\{\bo d_1^0,\bo d_2^0,\bo d_3^0\}$ (but the
$\bo d_i^0$ components are relative to $\{\bo e_1,\bo e_2,\bo
e_3\}$). Similarly, the $O(1)$ term of the angular momentum equation
(\ref{eq:amomentum}), projected on the director basis $\{\bo
d_1^0,\bo d_2^0,\bo d_3^0\}$ gives

{\setlength\arraycolsep{2pt}
\begin{eqnarray}
(M_1^0)'&=&\frac{2fM_3^0M_2^0}{\Gamma(1+R)}-\frac{fM_2^0M_3^0}{R}+F_2^0
+P\omega^2d_{23}^0(d_{22}^0d_{11}^0-d_{21}^0d_{12}^0),\label{stat:moments1}\\
(M_2^0)'&=&-\frac{2fM_3^0M_1^0}{\Gamma(1+R)}+fM_1^0M_3^0-F_1^0
+PR\omega^2d_{13}^0(d_{21}^0d_{12}^0-d_{11}^0d_{22}^0),\label{stat:moments2}\\
(M_3^0)'&=&\frac{fM_2^0M_1^0}{R}-fM_1^0M_2^0+PR\omega^2d_{13}^0
(d_{12}^0d_{31}^0-d_{11}^0d_{32}^0)+P\omega^2d_{23}^0(d_{22}^0d_{31}^0-d_{21}^0d_{32}^0).\label{stat:moments3}
\end{eqnarray}}

The $O(1)$ term of equation (\ref{eq:tangent}) can be expressed as
\begin{equation}
(\bo x^0)'=\bo d_{3}^0, \label{stat:tangent}
\end{equation}
and the twist equation (\ref{eq:twist2}) by
\begin{equation}
 (\bo d_{i}^0)'=\bo \kappa^0 \times \bo d_{i}^0, \qquad i=1,2,3,
\label{stat:twist}
\end{equation}
where $\bo \kappa^0 = \sum_{j=1}^3 \kappa_j^0 \bo d_{j}^0$. The
$O(1)$ term of constitutive relations (\ref{eq:constitutive}) can be
expressed as
\begin{equation}
M_1^0=\frac{1}{f}\kappa_1^0, \quad\quad
M_2^0=\frac{R}{f}\kappa_2^0, \quad\quad
M_1^0=\frac{\Gamma(1+R)}{2f}\kappa_3^0, \label{stat:constitutive}
\end{equation}
which can be used to express the $\kappa_i^0$ in (\ref{stat:twist})
in terms of the moments $M_i^0$.

Proceeding in the same way, the $O(1)$ part of the boundary conditions
is given by
{\setlength\arraycolsep{2pt}
\begin{eqnarray}
x^0 (0)= 0, \label{eq:et:stat:BC23}\\
y^0 (0)= 0, \label{eq:et:stat:BC24}\\
\bo F^0 (0) \cdot \bo e_3 = T, \label{eq:et:stat:BC25}\\
\bo d_3^0 (0) \cdot \bo e_1 = 0, \label{eq:et:stat:BC26}\\
\bo d_3^0 (0)\cdot \bo e_2 = 0, \label{eq:et:stat:BC27}\\
\bo d_1^0 (0) \cdot \bo e_2 = 0, \label{eq:et:stat:BC28}
\end{eqnarray}}
{\setlength\arraycolsep{2pt}
\begin{eqnarray}
x^0 (1)= 0, \label{eq:et:stat:BC29}\\
y^0 (1)= 0, \label{eq:et:stat:BC30}\\
z^0 (1)= 1, \label{eq:et:stat:BC31}\\
\bo d_3^0 (1) \cdot \bo e_1 = 0, \label{eq:et:stat:BC32}\\
\bo d_3^0 (1) \cdot \bo e_2 = 0, \label{eq:et:stat:BC33}\\
\bo d_1^0 (1) \cdot \bo e_2 = 0, \label{eq:et:stat:BC34}
\end{eqnarray}}
\begin{equation}
\begin{split}
\bo d_1^0 (0) \cdot \bo d_1^0 (0) = 1, \\
\bo d_2^0 (0) \cdot \bo d_2^0 (0) = 1, \\
\bo d_3^0 (0) \cdot \bo d_3^0 (0) = 1, \\
\bo d_1^0 (0) \cdot \bo d_2^0 (0) = 0, \\
\bo d_1^0 (0) \cdot \bo d_3^0 (0) = 0, \\
\bo d_2^0 (0) \cdot \bo d_3^0 (0) = 0, \\
\end{split}
\label{eq:et:stat:orthonorm}
\end{equation}

\subsection{The $O(\delta)$ equations -- linearisation}


The $O(\delta)$ part of the linear momentum equation (\ref{eq:lmomentum})
can be written as
\begin{equation}
\begin{split}
(\bo F^t(s,t))'+& \bo B_1(s)\bo F^t(s,t)+\bo B_2(s)\bo x^t(s,t)+\bo
B_3(s)\bo \alpha'(s,t)+\bo B_4(s)\bo \alpha(s,t)\\
=&\bo B_5(s) \ddot{\bo x}^t(s,t)+\bo B_6(s) \dot{\bo x}^t(s,t),
\end{split} \label{per:lmomentum}
\end{equation}
where the $3\times 3$ matrices $\bo B_i(s)$ are given in the
Appendix. Here we have expressed $\bo F^t$ relative to
$\{\bo d_1^0,\bo d_2^0,\bo d_3^0\}$ and $\bo x^t$ relative to
$\{\bo e_1,\bo e_2,\bo e_3\}$. For the $O(\delta)$ part of the
angular momentum equation (\ref{eq:amomentum}) we can write
\begin{equation}
\begin{split}
(\bo M^t(s,t))'+& \bo C_1(s)\bo M^t(s,t)+\bo C_2(s)\bo
\alpha'(s,t)+\bo C_3(s)\bo \alpha(s,t)+ \bo C_4(s)\bo F^t(s,t)\\
=& \bo C_5(s) \ddot{\bo \alpha}(s,t)+\bo C_6(s) \dot{\bo
\alpha}(s,t),
\end{split} \label{per:amomentum}
\end{equation}
where the matrices $\bo C_i(s)$ are again given in the Appendix.
$\bo M^t$ is expressed relative to $\{\bo d_1^0,\bo d_2^0,\bo d_3^0\}$.
The 9 twist equations (\ref{eq:twist2}) at $O(\delta)$ are reduced to
only 3  independent equations that relate $\bo \kappa^t$ and $\bo
\alpha$ as
\begin{equation}
\bo \kappa^t(s,t)=\bo \alpha'(s,t) +\bo \kappa^0(s) \times \bo
\alpha(s,t). \label{per:twist}
\end{equation}
Introducing these relations into the $O(\delta)$ part of the
constitutive relations gives
\begin{equation}
\bo M^t(s,t) + \bo D_1(s)\bo \alpha'(s,t)+\bo D_2(s)\bo \alpha(s,t)
= \bo D_3(s) \dot{\bo \alpha}(s,t)+\bo D_4(s)
\frac{\partial^2}{\partial s \partial t}(\bo \alpha(s,t)),
\label{per:constitutive}
\end{equation}
where the matrices $\bo D_i(s)$ are given in the Appendix. Finally, the
$O(\delta)$ part of equation (\ref{eq:tangent}) yields
\begin{equation}
\begin{split}
(\bo x^t(s,t))' = \alpha(s,t)\times\bo d_3^0(s).
\label{per:tangent}
\end{split}
\end{equation}

Applying the perturbation scheme to the boundary conditions at
$O(\delta)$, we obtain
{\setlength\arraycolsep{2pt}
\begin{eqnarray}
x^t(0,t)&=& 0,\label{per:et:bc1}\\
y^t(0,t)&=& 0,\label{per:et:bc2}\\
\bo F^t (0,t) \cdot \bo e_3 &=& 0, \label{eq:et:bc4}\\
\bo \alpha(0,t)&=&\bo 0,\label{per:et:bc3}\\
\bo x^t(1,t)&=& \bo 0,\label{per:et:bc5}\\
\bo \alpha(1,t)&=&\bo 0.\label{per:et:bc6}
\end{eqnarray}}

After elimination of the $\kappa_i^0$ by means of (\ref{stat:constitutive}),
the set of 12 equations (\ref{per:lmomentum}), (\ref{per:amomentum}),
(\ref{per:constitutive}) and (\ref{per:tangent}) together with the 12
boundary conditions (\ref{per:et:bc1})--(\ref{per:et:bc6}), with appropriate
initial conditions, form a well-posed initial-boundary-value problem.
%

\section{Stability analysis}
\label{sect:stab_anal}


Since we are interested in stability of solutions we look for solutions of
the $O(\delta)$ equations of the form
\begin{eqnarray}
\bo x^t (s,t) = \hat{\bo x}^t (s) e^{\lambda t},\label{per:sol1}\\
\bo \alpha (s,t) = \hat{\bo \alpha} (s) e^{\lambda t},\\
\bo F^t (s,t) = \hat{\bo F}^t (s) e^{\lambda t},\\
\bo M^t (s,t) = \hat{\bo M}^t (s) e^{\lambda t}.\label{per:sol4}
\end{eqnarray}
When these expressions are inserted into
(\ref{per:lmomentum})--(\ref{per:tangent}) a linear eigenvalue problem for
a 12-dimensional ODE is obtained in terms of the variables
$(\hat{\bo x}^t,\hat{\bo \alpha}^t,\hat{\bo F}^t,\hat{\bo M}^t)$.
The eigenvalue $\lambda$ measures the growth of small perturbations and is
to be found as past of the solution. Eigenvalues come as complex conjugate
pairs. A whirling state is unstable if at least one of the (in general
infinitely many) $\lambda$'s has positive real part.

To solve a real system of equations we split the eigenvalues and variables
(eigenfunctions) into real and imaginary parts,
$\lambda=\lambda_r + i \lambda_i$, $ \hat{\bo x}^t = \hat{\bo x}^t_r
+ i \hat{\bo x}^t_i$, $ \hat{\bo \alpha}^t = \hat{\bo \alpha}^t_r +
i \hat{\bo \alpha}^t_i$, $ \hat{\bo F}^t = \hat{\bo F}^t_r + i
\hat{\bo F}^t_i$ and $ \hat{\bo M}^t = \hat{\bo M}^t_r + i \hat{\bo
M}^t_i$. The equations (\ref{per:lmomentum})--(\ref{per:tangent}), along with
the boundary conditions (\ref{per:et:bc1})--(\ref{per:et:bc6}), are similarly
split into real and imaginary parts. Thus we end up with a doubled
24-dimensional linearised boundary-value problem.

\subsection{Stability of the straight rod -- magnetic buckling}
\label{sect:magn_buckling}

The trivial solution of the $O(1)$ equations
(\ref{stat:forces1})--(\ref{stat:constitutive}), representing a straight and
untwisted rod, is given by
\begin{equation}
\bo x(s)=s\bo e_3, \quad \bo F(s)=-T\bo e_3, \quad \bo M(s)=\bo
0, \quad \bo d_i(s)=\bo e_i~~~(i=1,2,3), \quad s\in [0,1].
\label{trivial_sol}
\end{equation}
For the statical case ($\omega=0$) without end force ($T=0$) the $O(\delta)$
equations (\ref{per:lmomentum}), (\ref{per:amomentum}),
(\ref{per:constitutive}), (\ref{per:tangent}) about this trivial solution,
on inserting (\ref{per:sol1})--(\ref{per:sol4}), can be written as
\begin{eqnarray}
&& x''''-\lambda^2Px''+\frac{\lambda^2}{R}x-\frac{B}{R}y'=0, \nonumber \\
&& y''''-\lambda^2Py''+\lambda^2y+Bx'=0, \label{lin_eqs} \\
&& M_3''-\frac{2\lambda^2P}{\Gamma}M_3=0, \nonumber
\end{eqnarray}
with boundary conditions
\begin{equation}
x(0)=x(1)=x'(0)=x'(1)=y(0)=y(1)=y'(0)=y'(1)=M_3'(0)=M_3'(1)=0,
\label{BCs}
\end{equation}
while $F_3\equiv 0$, $z\equiv 0$. Note that the torsional ($M_3$) modes
decouple from the bending ($x,y$) modes.

\begin{figure}
\begin{center}
\includegraphics[width=0.45\linewidth]{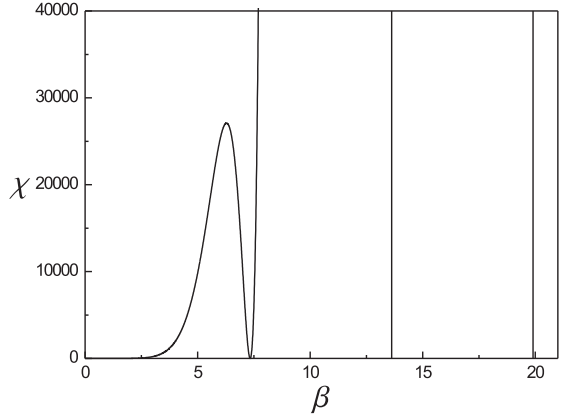}
\end{center}
\caption{Degenerate zeroes of the characteristic equation for the welded rod.
The first three critical $\beta$ values are 7.332130, 13.613561 and
19.896753.}
\label{fig:char_clamp}
\end{figure}

To find the statical magnetic buckling loads we set $\lambda=0$. The bending
equations then reduce to
\begin{equation}
z''''''+\frac{B^2}{R}z=0, \quad\quad \mbox{for} \quad\quad z=x',
\end{equation}
subject to
\begin{equation}
x(0)=x(1)=x'(0)=x'(1)=x''''(0)=x''''(1)=0.
\label{BC_x}
\end{equation}
On setting $z=e^{iks}$ we obtain the characteristic equation $-k^6+B^2/R=0$
with solutions
$k_{1,2}=\pm\beta,~k_{3,4,5,6}=\pm\beta\left(1\pm i\sqrt{3}\right)/2$,
where $\beta=B^{1/3}/R^{1/6}$. Application of the boundary conditions
(\ref{BC_x}) to the general solution $z(s)=\sum_{j=1}^6 a_je^{ik_js}$ leads
to the following critical condition:
\begin{eqnarray}
\chi(\beta)&:=&2\cos\beta + \cos 2\beta-2\left(\cos\frac{\beta}{2}+
\cos\frac{3\beta}{2}\right)\cosh\frac{\sqrt{3}\beta}{2}
+(2-\cos\beta)\cosh\sqrt{3}\beta \nonumber \\
&& -\sqrt{3}\sin\beta\,\sinh\sqrt{3}\beta-2\sqrt{3}\left(\sin\frac{\beta}{2}
-\sin\frac{3\beta}{2}\right)\sinh\frac{\sqrt{3}\beta}{2}=0.
\label{char_welded}
\end{eqnarray}
The critical loads correspond to pitchfork bifurcations where non-trivial
solutions bifurcate from the trivial straight solution. A plot of $\chi$
(Figure~\ref{fig:char_clamp}) shows that in the welded case the pitchfork
bifurcations are (doubly) degenerate, as was also found by Wolfe
\cite{wolfe3}. We stress that the above calculation is only possible for the
statics case. If $\omega\neq 0$ then the $x$ and $y$ equations do not
decouple and no simple characteristic equation is obtained.

\subsection{Eigenvalues for the unperturbed problem ($T=0$, $\gamma=0$,
$\omega=0$, $B=0$)}
\label{sect:unperturbed}

We shall call the case where $T=0$, $\gamma=0$, $\omega=0$ and $B=0$ the
{\it unperturbed problem}. For this problem explicit expressions can be
obtained for the eigenvalues of the linearisation about the straight solution.
The $x$ and $y$ equations in (\ref{lin_eqs}) decouple into two fourth-order
beam equations:
\begin{equation}
\begin{split}
x''''-\lambda^2Px''+\frac{\lambda^2}{R}x=0, \vspace*{0.2cm} \\
y''''-\lambda^2Py''+\lambda^2y=0,
\end{split}
\label{lin_eqs_0}
\end{equation}
subject to boundary conditions (\ref{BCs}). Since we anticipate imaginary
eigenvalues we set $\lambda=i\mu$, $x=e^{iks}$, $y=e^{i\kappa s}$, and find
for the $x$ equation
\begin{eqnarray}
&& k_{1,2}=\pm\left(\frac{1}{2}\mu^2P+\frac{1}{2}
\sqrt{\mu^4P^2+4\mu^2/R}\right)^{1/2}=:\pm a, \nonumber \\
&& k_{3,4}=\pm i\left(\frac{1}{2}\sqrt{\mu^4P^2+4\mu^2/R}-
\frac{1}{2}\mu^2P\right)^{1/2}=:\pm i b, \nonumber
\end{eqnarray}
while for the $y$ equation
\begin{eqnarray}
&& \kappa_{1,2}=\pm\left(\frac{1}{2}\mu^2P+\frac{1}{2}
\sqrt{\mu^4P^2+4\mu^2}\right)^{1/2}=:\pm\alpha, \nonumber \\
&& \kappa_{3,4}=\pm i\left(\frac{1}{2}\sqrt{\mu^4P^2+4\mu^2}-
\frac{1}{2}\mu^2P\right)^{1/2}=:\pm i\beta, \nonumber
\end{eqnarray}
where $a$, $b$, $\alpha$, $\beta$ are non-negative real numbers. The general
solutions are
\begin{eqnarray}
&& x(s)=A_x\sin as + B_x\cos as + C_x\sinh bs + D_x \cosh bs, \nonumber \\
&& y(s)=A_y\sin \alpha s + B_y\cos \alpha s + C_y\sinh \beta s +
D_y \cosh \beta s. \nonumber
\end{eqnarray}
Application of the boundary conditions (\ref{BCs}) leads to
\begin{eqnarray}
&& (a^2-b^2)\sin a\sinh b=2ab(1-\cos a\cosh b), \nonumber \\
&& (\alpha^2-\beta^2)\sin \alpha\sinh \beta=2\alpha\beta(1-\cos
\alpha \cosh \beta). \nonumber
\end{eqnarray}
In order to obtain the eigenvalues these two transcendental equations can
be solved numerically using a Newton-Raphson scheme. Meanwhile, the torsional
eigenvalues of the $M_3$ equation in (\ref{lin_eqs}) are given by
\[
\mu=\pm n\pi\sqrt{\frac{\Gamma}{2P}}.
\]
These are all the eigenvalues for the unperturbed problem. They will be used
as starting values in the numerical procedure described next.

\subsection{Numerical procedure}

The main idea is to use the known eigenvalues in the unperturbed problem as
starting values in a continuation procedure in order to compute the
eigenvalues and corresponding eigenfunctions for general values of the
parameters $T$, $\gamma$, $\omega$ and $B$. For this we use the well-tested
code AUTO \cite{doedel} (specifically AUTO2000). AUTO solves boundary-value
problems by means of orthogonal collocation. It requires a starting solution
and can then trace out solution curves as a parameter of the problem is
varied. Bifurcations are detected where branches of solutions intersect.
At such points AUTO is able to switch branches and compute curves of
bifurcating solutions.

Our procedure takes advantage of the fact that $\lambda$ appears only
quadratically in the linearisation (\ref{per:lmomentum}),
(\ref{per:amomentum}), (\ref{per:constitutive}) and (\ref{per:tangent}) if
$\gamma=0$ and $\omega=0$. To explain the method consider the typical
$O(\delta)$ equation
\begin{equation}
z''''-\lambda^2f(s)z''+\lambda^2g(s)z=0,
\end{equation}
where $f$ and $g$ are functions of the $O(1)$ solution. Writing $z=x+iy$,
$\lambda=\lambda_r+i\lambda_i$, we can decompose the $z$ equation into
\begin{equation}
\begin{split}
x''''-(\lambda_r^2-\lambda_i^2)f(s)x''+2\lambda_r\lambda_i f(s)y''+
(\lambda_r^2-\lambda_i^2)g(s)x-2\lambda_r\lambda_i g(s)y=0, \\
y''''-(\lambda_r^2-\lambda_i^2)f(s)y''-2\lambda_r\lambda_i f(s)x''+
(\lambda_r^2-\lambda_i^2)g(s)y+2\lambda_r\lambda_i g(s)x=0.
\end{split}
\label{re_im_eq}
\end{equation}
The important thing to note here is that these equations decouple into two
identical equations if the eigenvalue is either imaginary ($\lambda_r=0$) or
real ($\lambda_i=0$).

This suggests the following sequence of steps, involving boundary-value
problems of increasing dimension, to compute eigenvalues of statical or
uniformly whirling solutions.

\begin{enumerate}

\item
Consider the unperturbed problem of Section~\ref{sect:unperturbed} and,
noting that all eigenvalues are purely imaginary, solve the 30-dimensional
system of 18 $O(1)$ equations and one 12-dimensional system for the
imaginary part of the $O(\delta)$ equations (cf. the $y$ equation in
(\ref{re_im_eq})). Set $\lambda_r=0$ and use $\lambda_i$ as the continuation
parameter in AUTO in order to compute the eigenvalues (instead of solving the
transcendental equations in Section~\ref{sect:unperturbed}). These eigenvalues
will show up as branching points (BP), or pitchfork bifurcations, as
eigenvalues by definition are those values for which non-zero BVP solutions
exist. By symmetry it is only necessary to consider $\lambda_i>0$.

\item
Keeping the same 30-dimensional system, switch branches at a BP to
compute (`grow') the corresponding (imaginary) eigenfunction. Since the
equations are linear the value of $\lambda_i$ will not change in this run.
For later use we monitor the non-zero solution by means of some measure
$||.||_i$ (not necesarily a proper norm) on the space
$\{\hat{\bo x}_i^t,\hat{\bo \alpha}_i^t,\hat{\bo F}_i^t,\hat{\bo M}_i^t\}$
of imaginary linearised variables.

\item
Now consider the full system of 42 equations (18 $O(1)$ equations and two
sets of 12-dimensional $O(\delta)$ equations (cf.~(\ref{re_im_eq})). Fix the
measure $||.||_i$ on the imaginary part and release $\lambda_r$ instead in
order to compute the real eigenfunction (since the imaginary part of the
solution is fixed there is only one branch of solutions through the starting
point and there is nowhere else to go for the continuation but to `grow' the
real eigenfunction). Again we monitor this function by means of a suitable
measure $||.||_r$. In this run neither $\lambda_r$ nor $\lambda_i$ will
change.

This approach works because the solution obtained in step 2 also solves
the full 42-dimensional system when the extra 12 variables
$(\hat{\bo x}_r^t,\hat{\bo \alpha}_r^t,\hat{\bo F}_r^t,\hat{\bo M}_r^t)$
are set to zero. This is a consequence of the fact that the real and imaginary
parts of the $O(\delta)$ equations decouple if $\lambda_r=0$, as a result
of the quadratic dependence of the eigenvalue problem on $\lambda$
(cf.~(\ref{re_im_eq})).

\end{enumerate}

Steps 2 and 3 can be performed for as many of the BPs computed in step 1 as
required and will give the corrsponding eigenvalues and eigenfunctions. Once
these have been obtained both measures $||.||_r$ and $||.||_i$ can be fixed
and an extra system parameter such as $B$ or $\omega$ released in order to
trace the eigenvalues (and hence monitor stability changes) as system
parameters are varied. (Note that fixing $||.||_r$ and $||.||_i$ makes sense
as eigenfunctions are only defined up to a multiplicative factor.)

The above 3-step procedure is not limited to linearisations about the trivial
straight solution. It can be applied to any {\it starting solution} that has
no eigenvalue with both $\lambda_r$ and $\lambda_i$ non-zero, as these would
not be picked up in step 1. (It is of course no problem if eigenvalues
become fully complex (e.g., in a Hopf bifurcation) in the course of further
continuations.) For instance, we find that at the first critical $B$, given
by (\ref{char_welded}), the lowest conjugate pair of eigenvalues
$\pm\lambda_i$ goes to zero and becomes a real pair of eigenvalues,
signalling a stability change of the straight rod. The (first-mode) solution
bifurcating at this point is stable with all eigenvalues being imaginary and
the above procedure can be applied to find the eigenvalues.

We end this section with a few comments:

\begin{enumerate}

\item[$(i)$]
There are infinitely many eigenvalues and the above procedure only finds the
lowest order ones. This is of course a limitation of any numerical scheme.
We find that eigenvalues vary slowly with system parameters, suggesting
that stability is governed by the lowest-order eigenvalues. We typically
consider 5 or 6 eigenvalues.

\item[$(ii)$]
Note that in steps 1 and 2 above we could not have taken the full
42-dimensional system of equations as that would have made the branching
points (pitchfork bifurcations) degenerate and AUTO would not detect a BP.
This is because if $\lambda_r=0$ (or $\lambda_i=0$) the two sets of
12-dimensional linearised equations are identical (cf.~(\ref{re_im_eq})).

\end{enumerate}

\section{Numerical results}\label{results}

\subsection{The statics case, $\omega=0$}

\begin{figure}
\begin{center}
\includegraphics[width=0.7\linewidth]{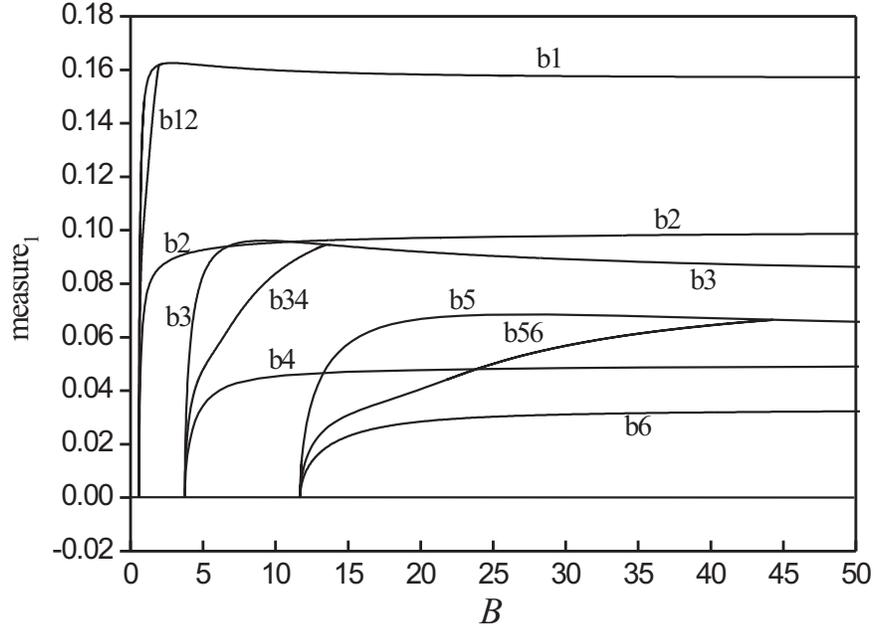}
\end{center}
\caption{Bifurcation diagram for the statical case ($\omega = 0$).
System parameters are those of Table~\ref{table:et:parameters}.
measure$_1=\int_0^{1}|x^0(s)|\textrm{d}s$.}
\label{fig:tet:bif-measure1}
\end{figure}

\begin{figure}
\begin{center}
\includegraphics[width=0.7\linewidth]{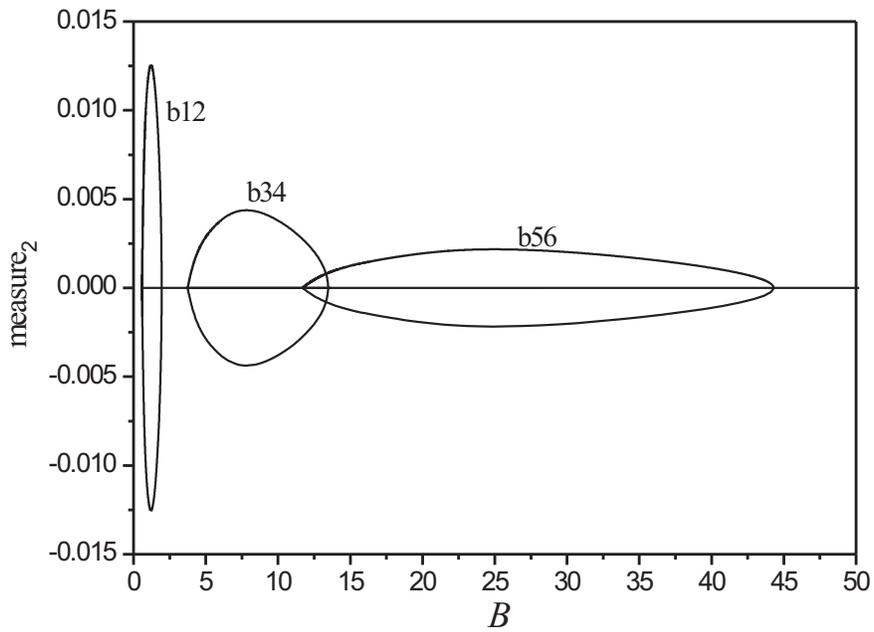}
\end{center}
\caption{Bifurcation diagram for the statical case ($\omega = 0$) using a
different solution measure. measure$_2=\int_0^{1}x^0(s)y^0(s)\textrm{d}s$.
The primary branches have measure$_2=0$.}
\label{fig:tet:bif-measure2}
\end{figure}

\begin{figure}
\begin{center}
\includegraphics[width=1\linewidth]{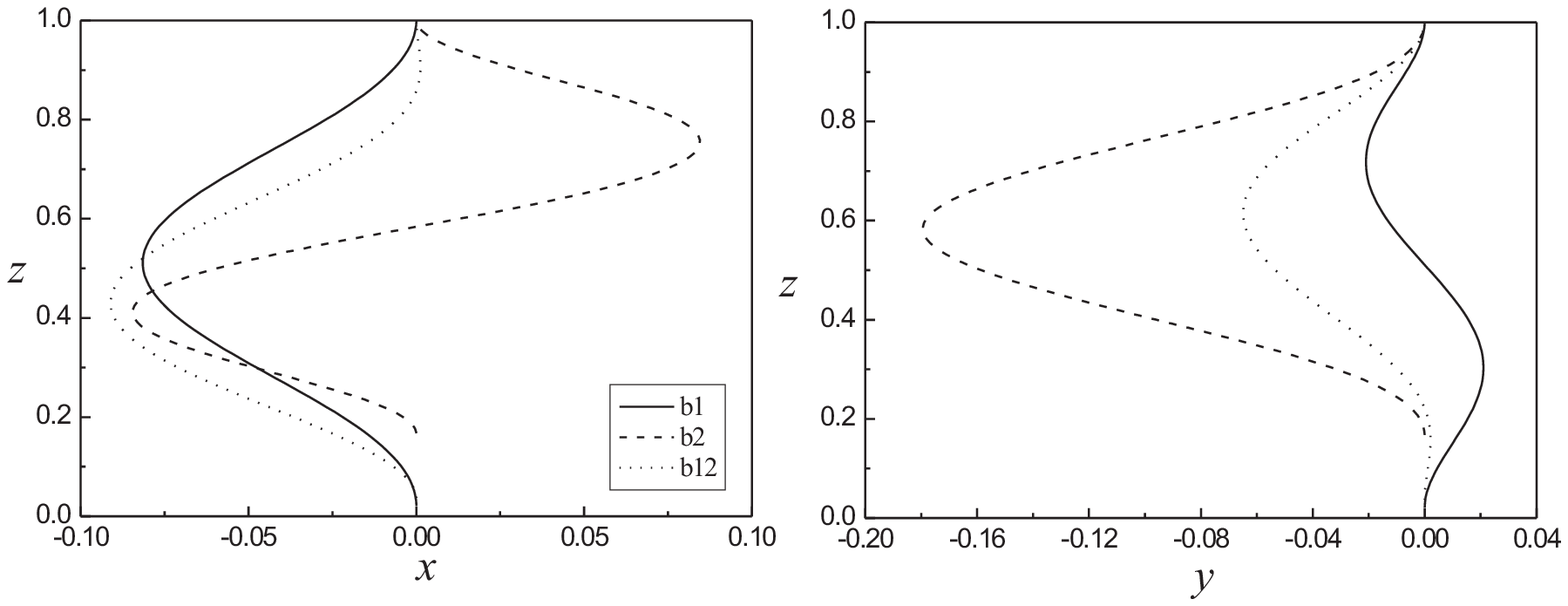}
\end{center}
\caption{Projection of the modes corresponding to branches b1, b2
and transition branch b12 onto the \{$\bo e_1$-$\bo e_3$\} and
\{$\bo e_2$-$\bo e_3$\} planes at constant measure$_1 = 0.04$.
Values of $B$ are: 0.595 (b1), 0.605 (b12) and 0.7 (b2).}
\label{fig:tet:tether-modes-transition}
\end{figure}

Figures~\ref{fig:tet:bif-measure1}~and~\ref{fig:tet:bif-measure2}
show the bifurcation diagram obtained when the magnetic field
parameter $B$ is varied, using different measures along the vertical
axes. The dimensionless parameters taken are those of the SET and
are listed in Table~\ref{table:et:parameters}, along with the
dimensional parameters. We observe that at critical values of $B$
pairs of non-trivial solution branches bifurcate, thus confirming
the degenerate nature of the bifurcation found analytically in
Section \ref{sect:magn_buckling}. The first three bifurcations occur
at $B=0.585\,(BP1)$, $3.747\,(BP2)$ and $11.698\,(BP3)$, in agreement with
Figure~\ref{fig:char_clamp} when the numerical factor $f$ is taken
into account. In addition to these primary branches, there are
branches, labelled b12, b34 and b56, connecting these pairs. One end
of these connecting branches reaches down to the horizontal axis,
while the other end comes out of a secondary pitchfork bifurcation
along the first branch of the pair of primary branches. They in fact
come in symmetric pairs as is brought out by the measure plotted in
Figure~\ref{fig:tet:bif-measure2}. The shapes of the rod along the
connecting branches form a smooth transition between the shapes on
the connected branches; see
Figure~\ref{fig:tet:tether-modes-transition} where projections onto
the \{$\bo e_1$-$\bo e_3$\} and \{$\bo e_2$-$\bo e_3$\} planes of
bifurcating solutions along branches b1, b2 and b12 are shown at
constant measure$_1 = 0.04$.

\begin{table}
  \centering
  \caption{Dimensional and dimensionless parameters for the SET.}
  \vspace{0.05\linewidth}
  \begin{tabular}{|c|c||c|c|} \hline
    $L$ & 100 m & $P$ & 0.001\\ \hline       
    $A$ & $2.879 \times 10^{-11}$ m$^2$ & $R$ & 0.5526\\ \hline
    $E$ & $1.32 \times 10^{11}$ N/m$^2$ & $\Gamma$ & 0.76923\\ \hline
    $E I_1$ & 38 Nm$^2$ & $f$ & 500.5639\\ \hline
    $E I_2$ & 21 Nm$^2$ & &\\ \hline
  \end{tabular}
  \label{table:et:parameters}
\end{table}


\begin{figure}
\begin{center}
\includegraphics[width=0.55\linewidth]{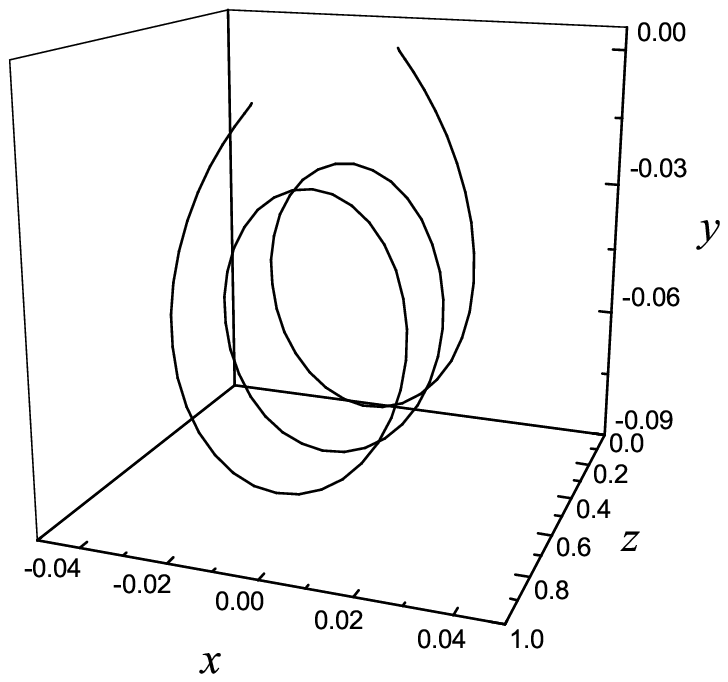}
\end{center}
\caption{3D view of a solution along branch b5.}
\label{fig:tet:b5-3D}
\end{figure}

Solutions bifurcating at larger $B$ values have successively more coils.
Figure~\ref{fig:tet:b5-3D} shows the three-dimensional shape
of a solution on the fifth bifurcating branch.


\begin{figure}
\begin{center}
{\bf~~~~~~~~~~(a)~~~~~~~~~~~~~~~~~~~~~~~~~~~~~~~~~~~~~~~~~~~~~~~~~(b)~~~}\\
\includegraphics[width=0.45\linewidth]{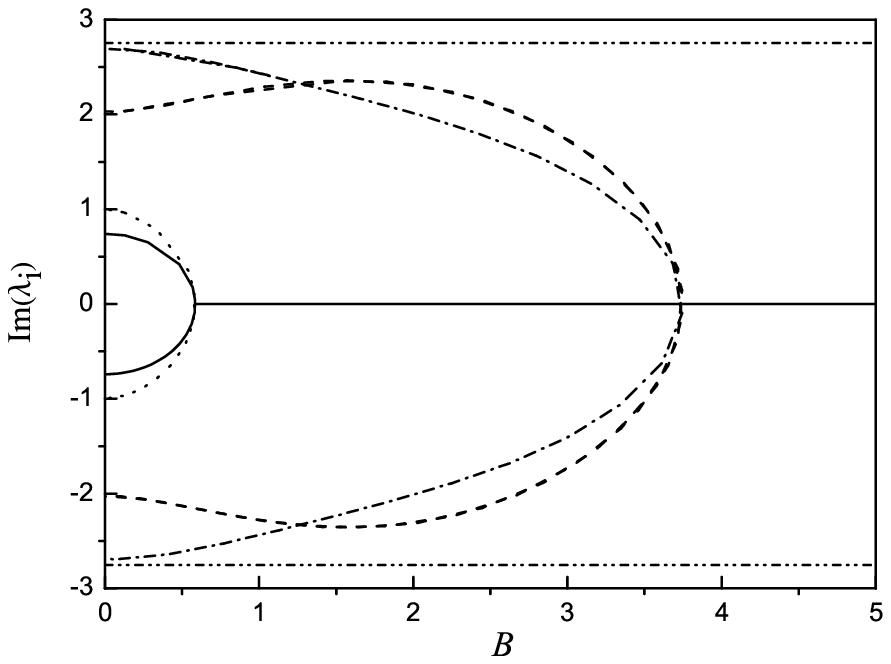}
~~~~~\includegraphics[width=0.45\linewidth]{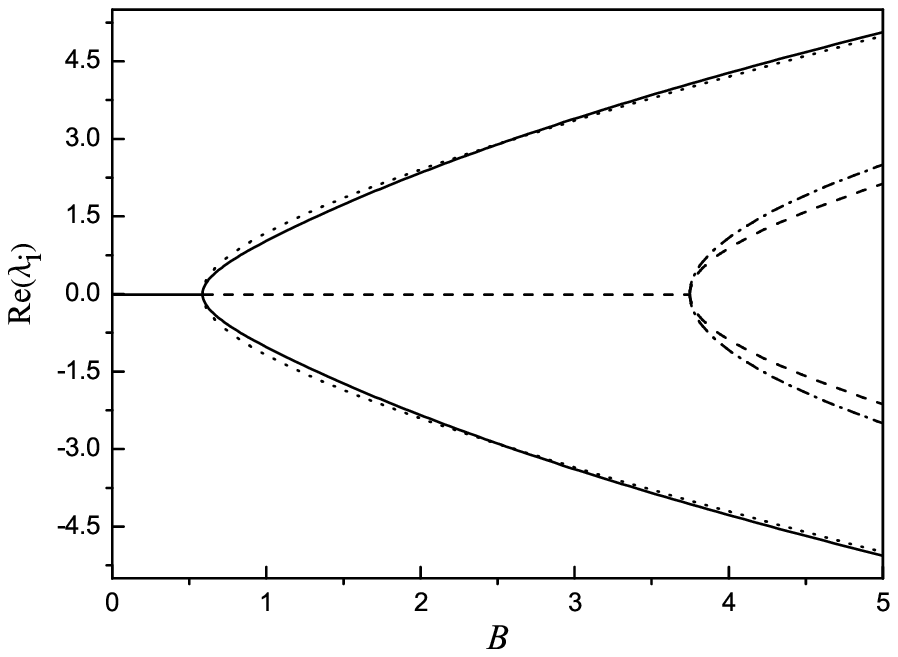}
\end{center}
\caption{Evolution of imaginary (a) and real (b) part of the first five
pairs of eigenvalues along the trivial solution as $B$ is increased from 0
to $>BP2$.}
\label{fig:tet:welded_w0-stab-0BP1BP2-imag}
\end{figure}

Figure~\ref{fig:tet:welded_w0-stab-0BP1BP2-imag} shows the evolution
of the imaginary and real parts of the first five pairs of
eigenvalues along the trivial solution, from $B=0$ to $B>BP2$ (i.e.,
beyond the second pitchfork). As the branches are born in pairs at
the degenerate pitchfork, two pairs of pure imaginary eigenvalues
collide at zero and become pure real pairs at $B=0.585\,(BP1)$. One
of the eigenvalues of each real pair is positive, meaning that the
trivial solution becomes unstable at $B=BP1$, as one would expect. A
similar further loss of stability occurs at $B=BP2$.

\begin{figure}
\begin{center}
\includegraphics[width=0.6\linewidth]{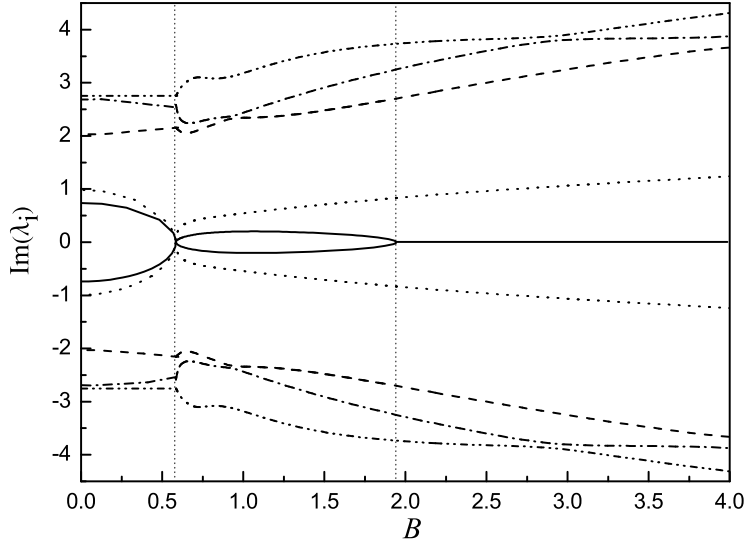}
\end{center}
\caption{Imaginary part evolution of the first five pairs of
eigenvalues along the trivial solution from $B=0$ to $B=BP1$, then
following branch b1 beyond the secondary pitchfork.}
\label{fig:tet:welded-eigen-0-BP1-b1-b12-imag}
\end{figure}

\begin{figure}
\begin{center}
\includegraphics[width=0.6\linewidth]{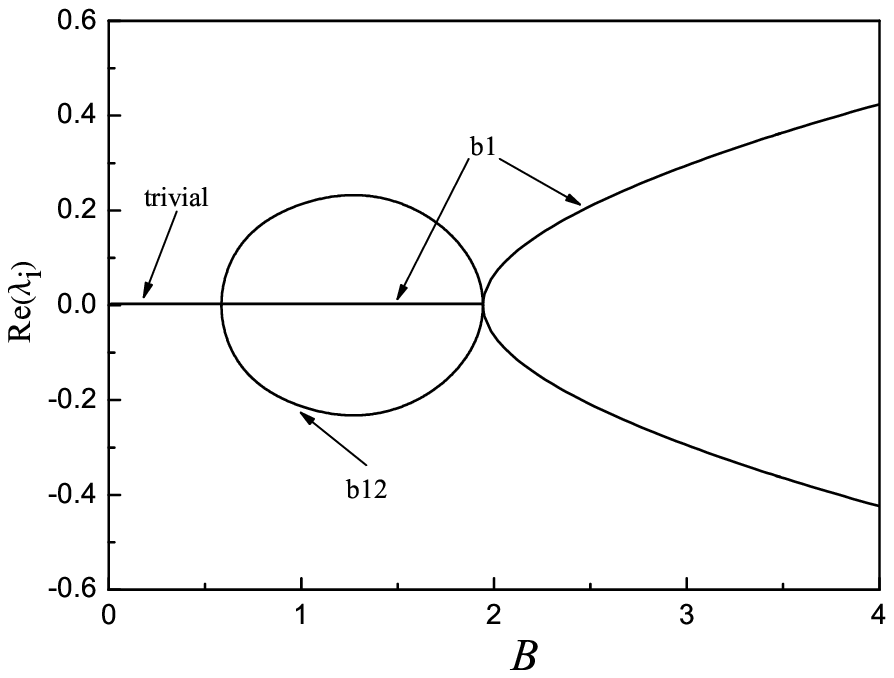}
\end{center}
\caption{Real part evolution of the first pair of eigenvalues along
the trivial solution from $B=0$ to $B=BP1$, then following branch b1
beyond the secondary pitchfork as well as connecting branch b12.}
\label{fig:tet:welded-eigen-0-BP1-b1-b12-real}
\end{figure}

Figure~\ref{fig:tet:welded-eigen-0-BP1-b1-b12-imag} shows the evolution of
the first five pairs of purely imaginary eigenvalues when switching
from the stable trivial branch to branch b1, which is initially found to be
stable. The figure reveals that branch b1 loses stability at the secondary
pitchfork bifurcation at $B=1.942$, where b12 connects. At this point the
first pair of eigenvalues becomes real, as illustrated in
Figure~\ref{fig:tet:welded-eigen-0-BP1-b1-b12-real}. This figure also shows
that connecting branch b12 is unstable as it too has a pair of real
eigenvalues.



Considering the problem of the ET from a design point of view, it is
interesting to quantify the value of the geomagnetic field at which the
straight tether would buckle. In the case at hand, this value is $B = 0.587$
(see Figure~\ref{fig:tet:bif-measure1}), which in dimensional parameters
yields $IB_0 = 2.104\times 10^{-5}~\frac{N}{m}$. Noting that the
maximum value of the geomagnetic field is $B_g = 7\times
10^{-5}~T$~\cite{lorenzini}, and assuming that the maximum current
that will flow along the tether would be $I = 1~A$~\cite{lorenzini},
the maximum expected value for the constant $IB_g = 7\times
10^{-5}~\frac{N}{m}$. This means that the ET would be below the
critical value at any time. Note, though, that $B$ goes as the cube
of $L$ and for longer tethers, which are common in radially
stabilised ETs, the critical value may be exceeded resulting in
buckling into a coiled shape (cf.~Figure~\ref{fig:tet:b5-3D}),
as has been reported in some tether flights such as the PGM and
TSS-1R missions~\cite{beletsky,lorenzini}.

\subsection{Whirling solutions (relative equilibria), $\omega\neq0$}

\begin{figure}
\begin{center}
\includegraphics[width=0.9\linewidth]{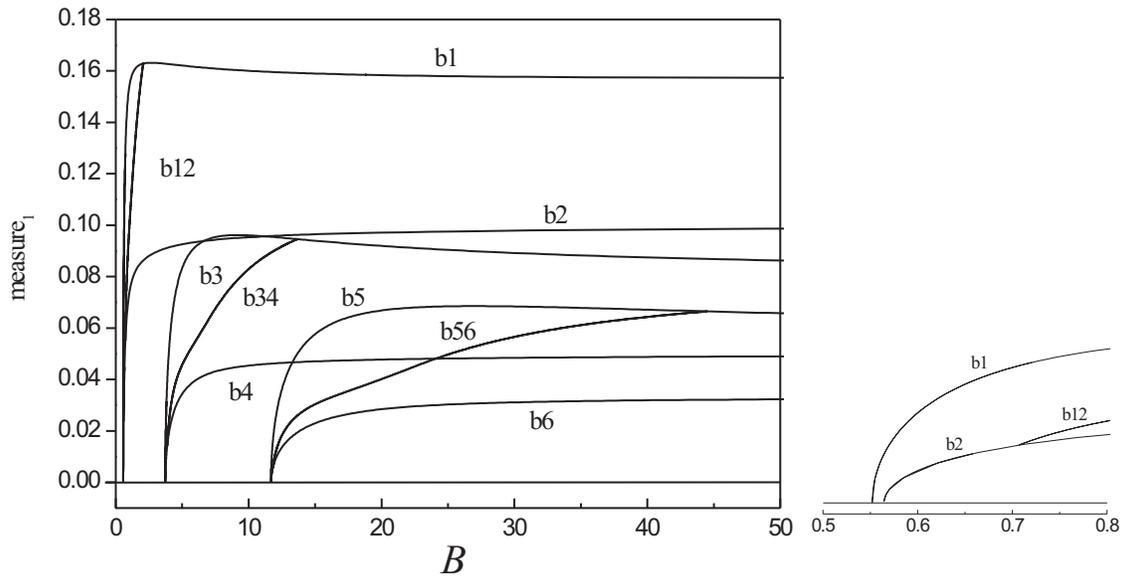}
\end{center}
\caption{Bifurcation diagram for a whirling anisotropic rod at $\omega = 0.25$.
Branches b1 and b2 still bifurcate (at distinct values of $B$) and trivial
solutions before the first pitchfork are stable.}
\label{fig:tet:bif-w025-measure1}
\end{figure}

\begin{figure}
\begin{center}
\includegraphics[width=0.9\linewidth]{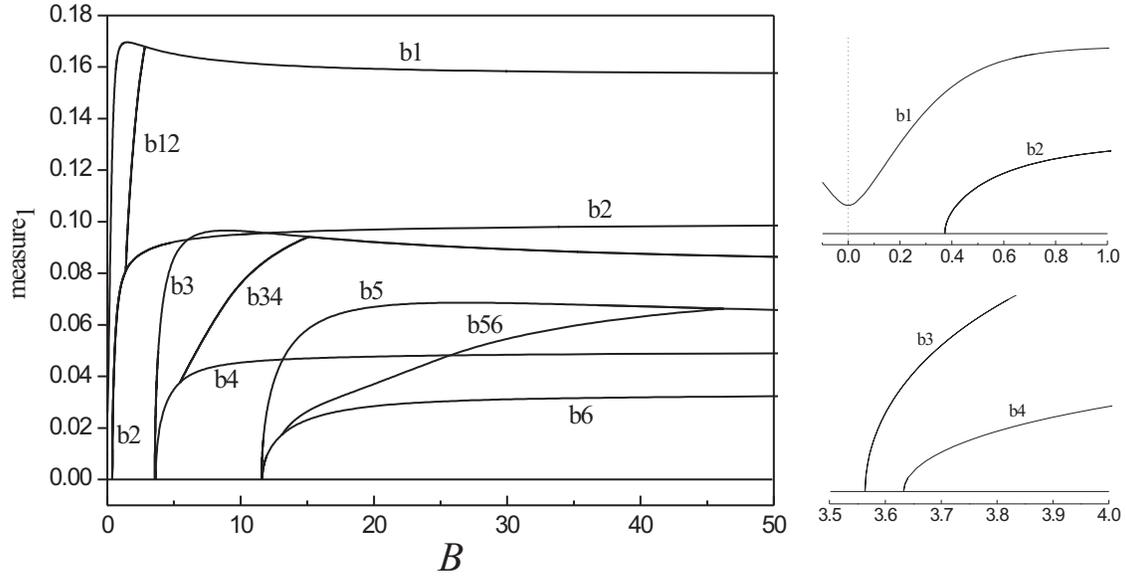}
\end{center}
\caption{Bifurcation diagram for a whirling anisotropic rod at $\omega = 0.75$.
Branch b1 no longer bifurcates and all trivial solutions are unstable.}
\label{fig:tet:bif-w075-measure1}
\end{figure}

\begin{figure}
\begin{center}
\includegraphics[width=0.7\linewidth]{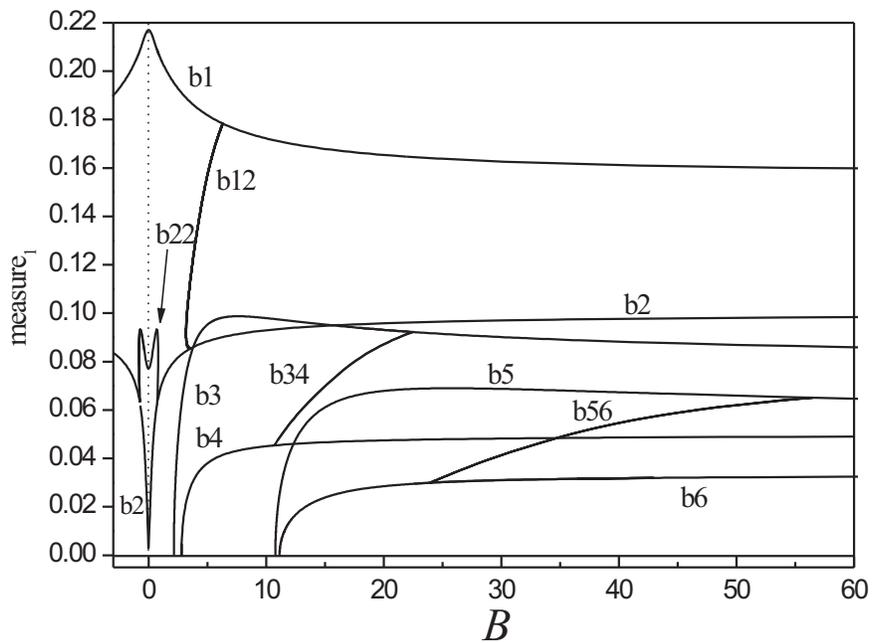}
\end{center}
\caption{Bifurcation diagram for a whirling anisotropic rod at $\omega = 2$.
Branches b1 and b2 no longer bifurcate. An extra (unstable) branch b22 has
appeared.}
\label{fig:tet:bif-w2-measure1}
\end{figure}

We now introduce angular velocity $\omega$ to the ET and seek
relative equilibria, that is, solutions that appear static when
viewed from the rotating frame \{$\bo e_1,\bo e_2,\bo e_3$\}.
Figures~\ref{fig:tet:bif-w025-measure1},
\ref{fig:tet:bif-w075-measure1} and \ref{fig:tet:bif-w2-measure1}
show the effect of increasing angular velocity on the bifurcation
diagram of Figure~\ref{fig:tet:bif-measure1}. The degenerate
pitchfork bifurcations get resolved and the connecting branches move
up along the principal branches. For $\omega=0.25$
(Figure~\ref{fig:tet:bif-w025-measure1}) all branches, b1, b2, ...,
bifurcate, as in the statics case of
Figure~\ref{fig:tet:bif-measure1}, and the trivial solution is
stable up to the first primary pitchfork. However, at $\omega=0.75$
(Figure~\ref{fig:tet:bif-w075-measure1}), b1 has merged with its
symmetric partner at $B<0$ and lifted off the horizontal axis,
leaving no stable straight solutions. b2 similarly lifts off near
$\omega=1.25$. Figure~\ref{fig:tet:bif-w2-measure1} shows the
bifurcation diagram at $\omega=2$. An extra branch of connecting
solutions (b22) has appeared, but no further branches have become
non-bifurcating.

\begin{figure}
\begin{center}
\includegraphics[width=0.7\linewidth]{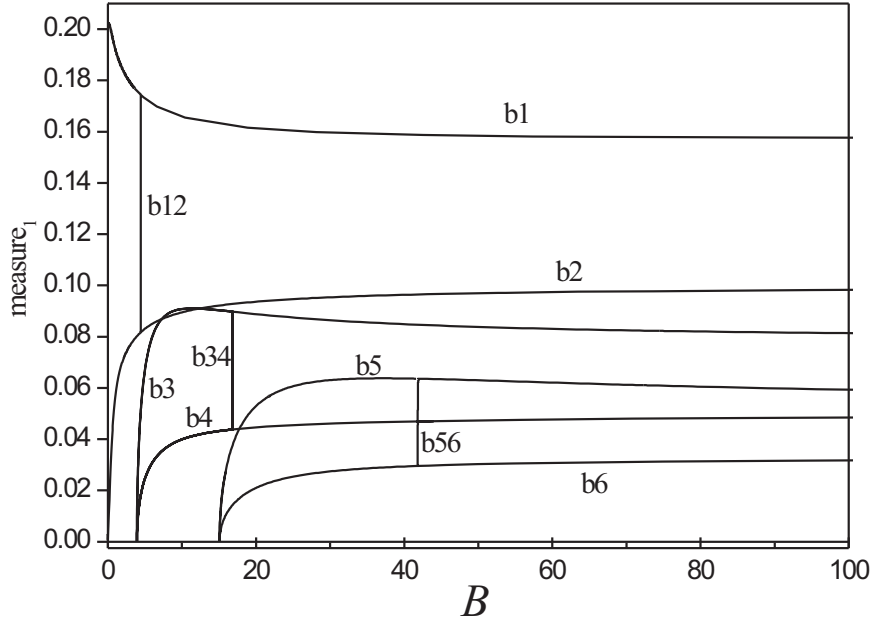}
\end{center}
\caption{Bifurcation diagram for a whirling isotropic rod ($\omega=2$,
$R=1$).}
\label{fig:tet:bif-w2-r1-measure2}
\end{figure}

Figure~\ref{fig:tet:bif-w2-r1-measure2} shows the bifurcation diagram at
$\omega=2$ for the case of an isotropic rod ($R=1$). Here the bifurcations
become degenerate again but the connecting branches do not move down. In
fact, they turn into vertical branches. We conclude that neither anisotropy
($R$) nor whirl ($\omega$) alone resolves the pitchfork degeneracy; both are
required.

\begin{figure}
\begin{center}
{\bf~~~~~~~~(a)~~~~~~~~~~~~~~~~~~~~~~~~~~~~~~~~~~~~~~~~~~~~~~~~~~~(b)~~~}\\
\includegraphics[width=0.45\linewidth]{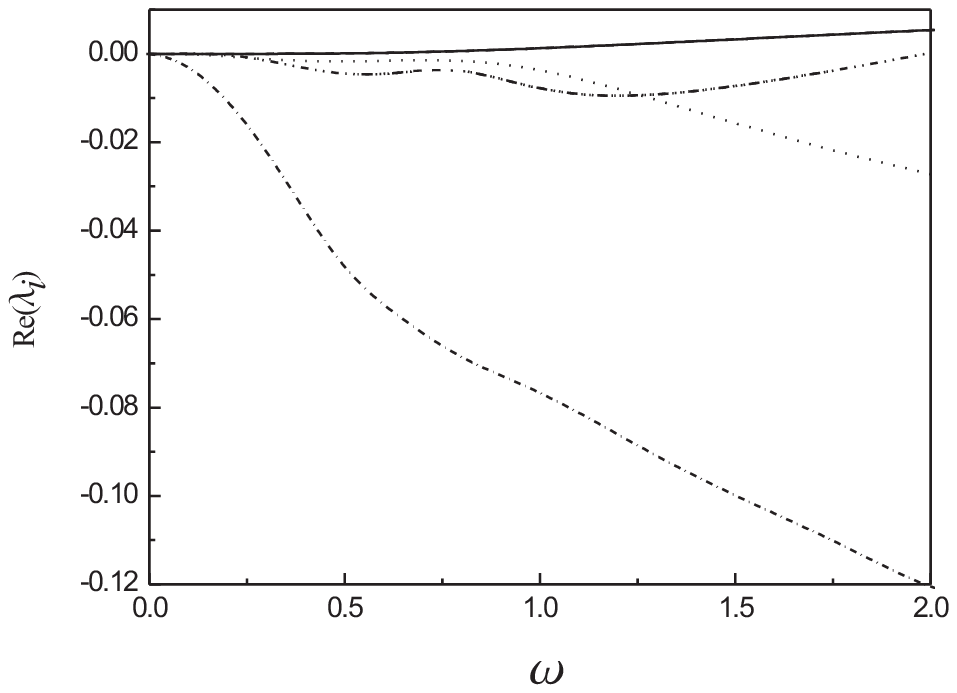}
~~~~~\includegraphics[width=0.45\linewidth]{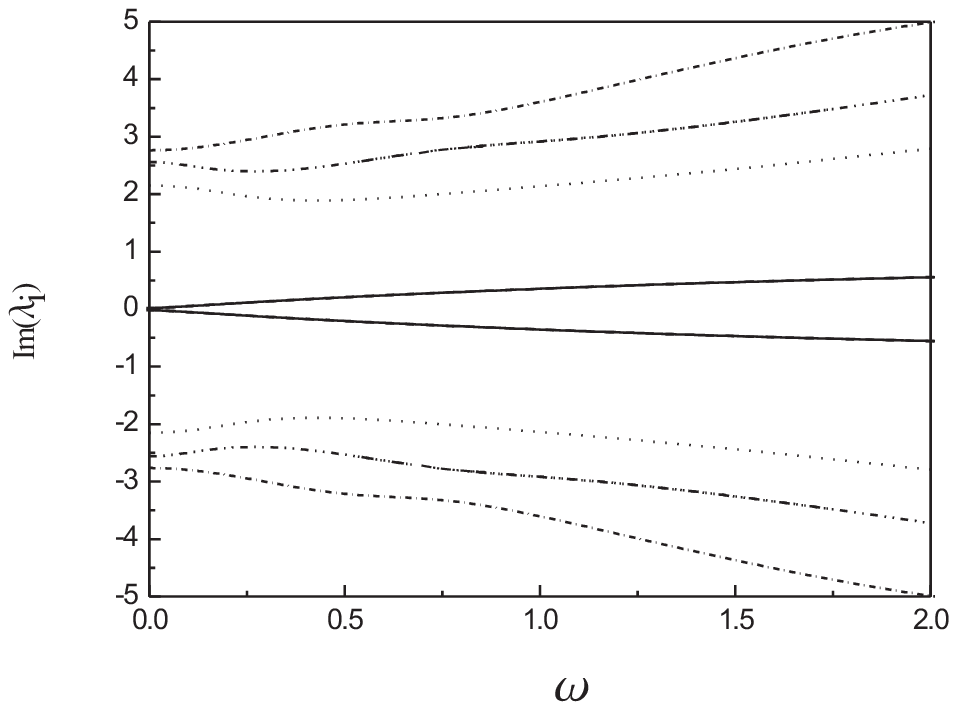}
\end{center}
\caption{Evolution with $\omega$ of real (a) and imaginary (b) parts of the
first five pairs of eigenvalues from branch b1 of
Figure~\ref{fig:tet:bif-measure1} ($\omega=0$) with $B=0.5869$.}
\label{fig:tet:eigenb1_B0-587_wcont-real}
\end{figure}


To investigate the stability of non-trivial solutions we follow the
eigenvalues computed earlier as $\omega$ is increased from 0.
Figure~\ref{fig:tet:eigenb1_B0-587_wcont-real} shows the evolution
of the real and imaginary parts of the first five conjugate pairs of
eigenvalues when continuing in $\omega$ from a static solution on the first
branch (b1) of Figure~\ref{fig:tet:bif-measure1} at fixed $B=0.5869$.
Note that the first and second eigenvalues are almost coincident and can
hardly be distinguished in the figure.
Both eigenvalues acquire positive real parts; hence the static coiled
tether becomes unstable as soon as it is being rotated about $\bo k$.
The other eigenvalues move to the negative real half-plane as $\omega$
is increased. The angular velocity is increased to $\omega=2$, which is
beyond the critical velocity in the $B=0$-case where the tether
buckles at the first bending natural frequency of the beam,
$\omega=\sqrt{R}$~\cite{valverde-cosserat}.

%

\begin{figure}
\begin{center}
\includegraphics[width=0.45\linewidth]{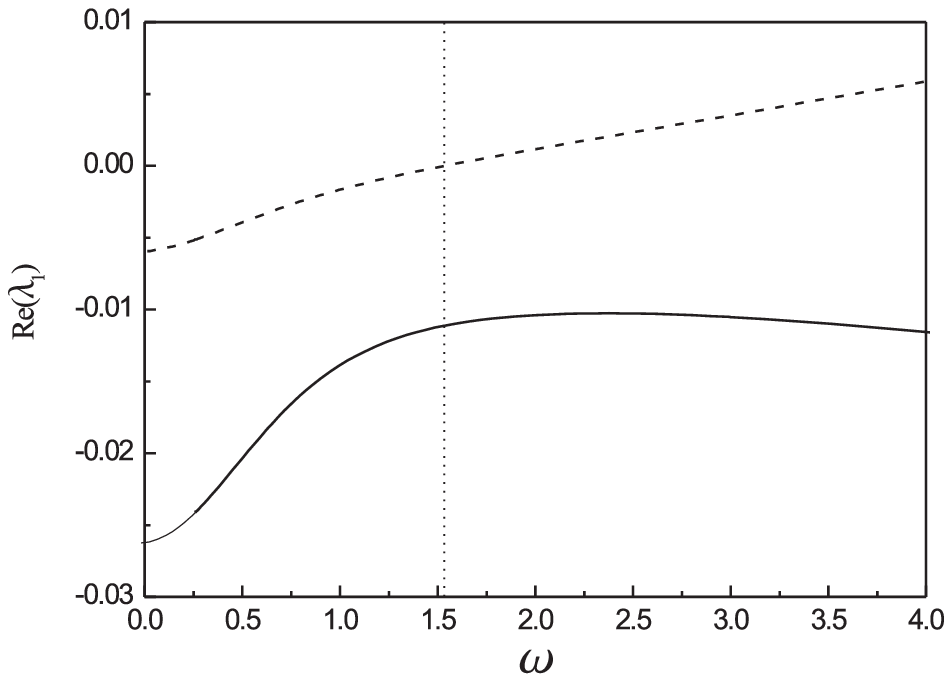}
~~~~~\includegraphics[width=0.45\linewidth]{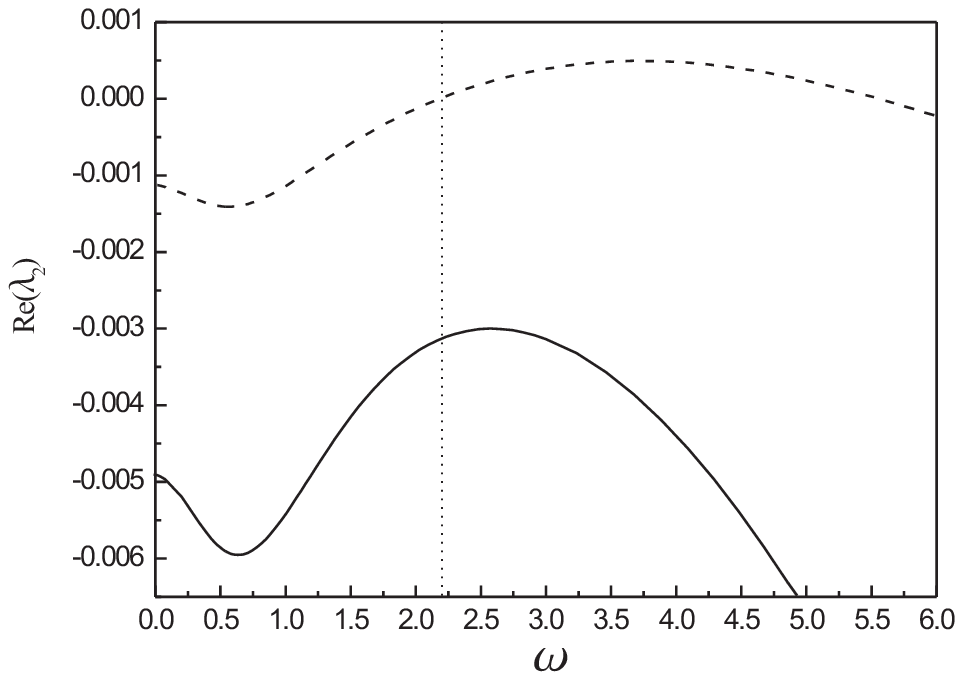} \vspace{0.5cm}\\
\includegraphics[width=0.45\linewidth]{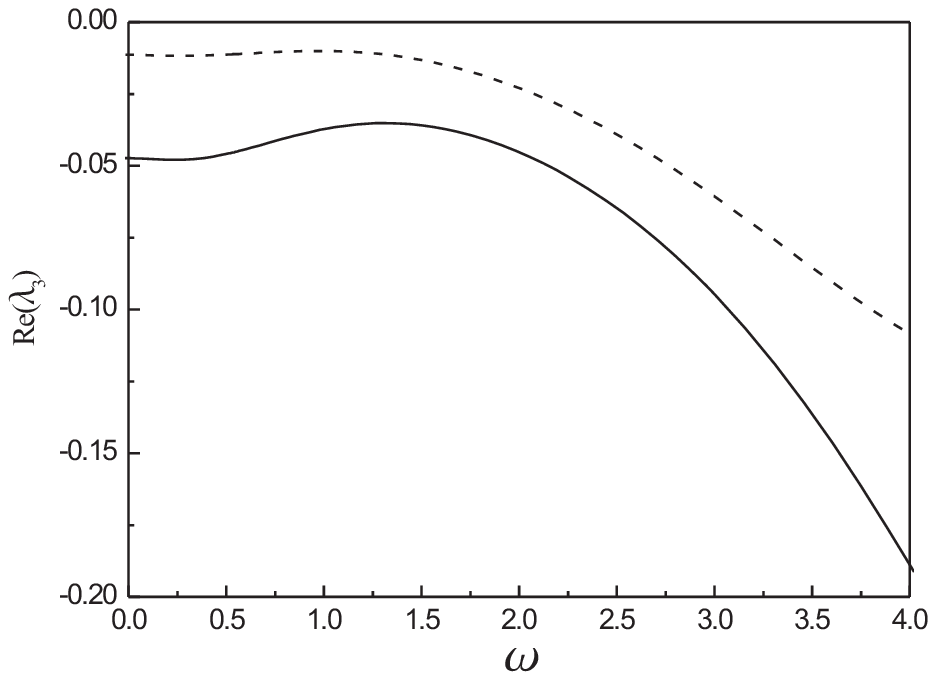}
\end{center}
\caption{Evolution with respect to $\omega$ of the real part of
the first three pairs of eigenvalues for a solution along b1 at
$\gamma = 0.04375$ (solid) and $\gamma = 0.01$ (dashed). ($B=1.2$.)}
\label{fig:tet:helix_stab-b1_B35_gamma05_w-real}
\end{figure}

\begin{figure}
\begin{center}
\includegraphics[width=0.45\linewidth]{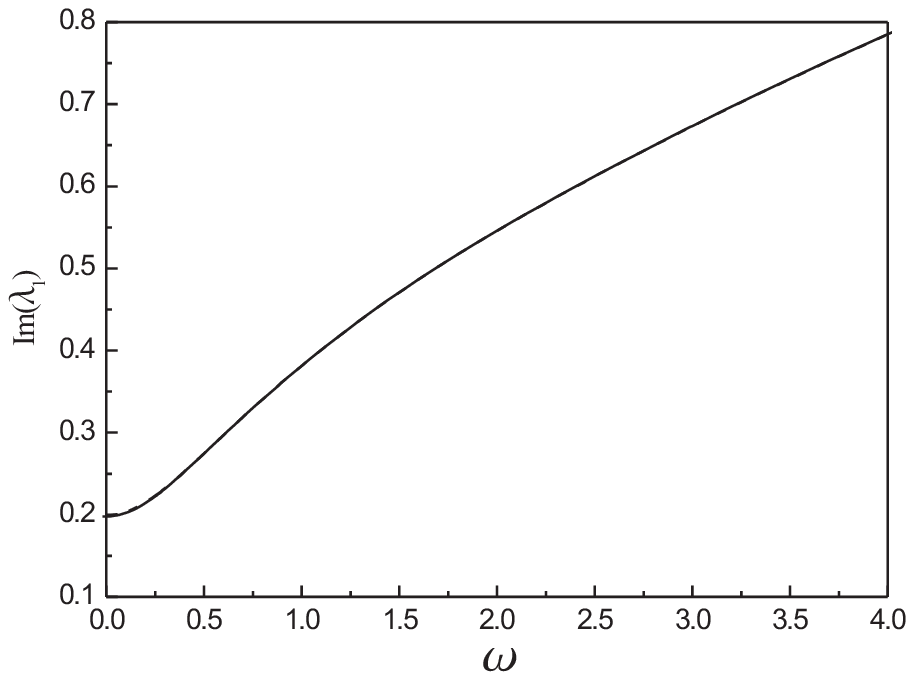}
~~~~~\includegraphics[width=0.45\linewidth]{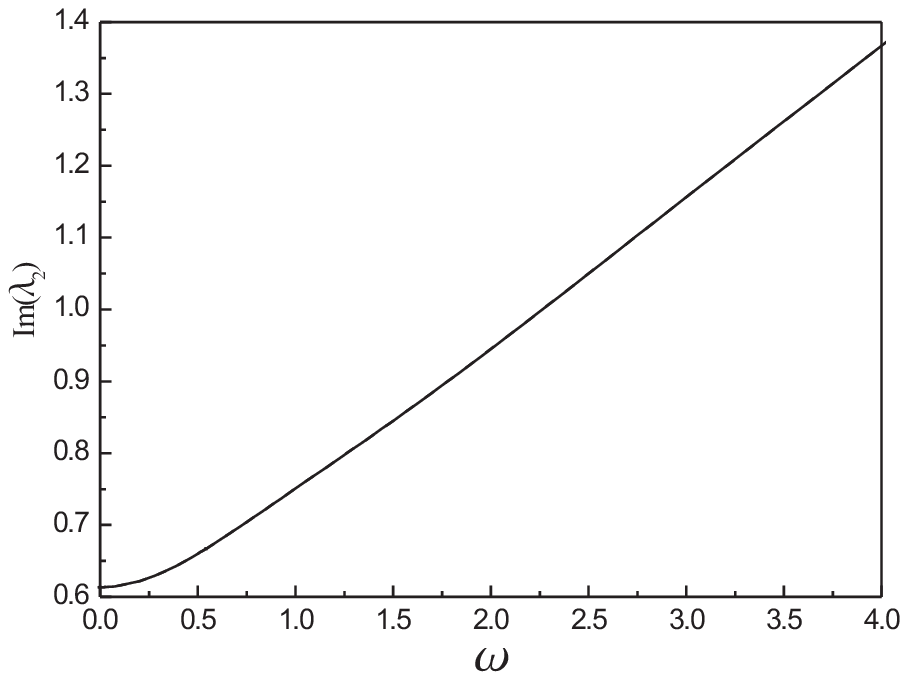} \vspace{0.2cm}\\
\includegraphics[width=0.45\linewidth]{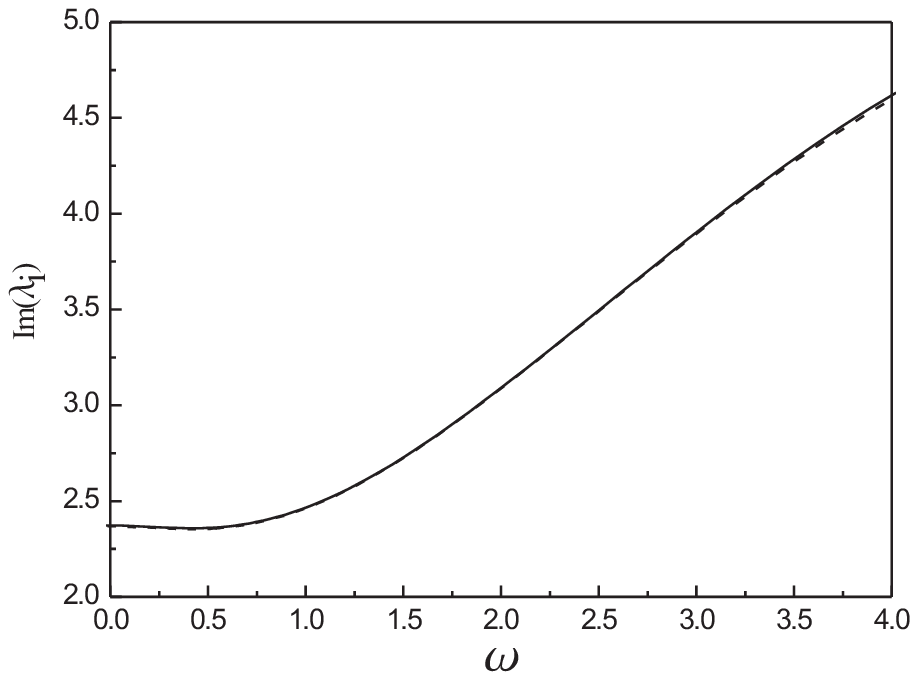}
\end{center}
\caption{Evolution with respect to $\omega$ of the imaginary part of
the first three pairs of eigenvalues for a solution along b1 at
$\gamma = 0.04375$ (solid) and $\gamma = 0.01$ (dashed). The curves are
independent of $\gamma$. ($B=1.2$.)}
\label{fig:tet:helix_stab-b1_B35_gamma05_w-imag}
\end{figure}


\begin{figure}
\begin{center}
\includegraphics[width=0.5\linewidth]{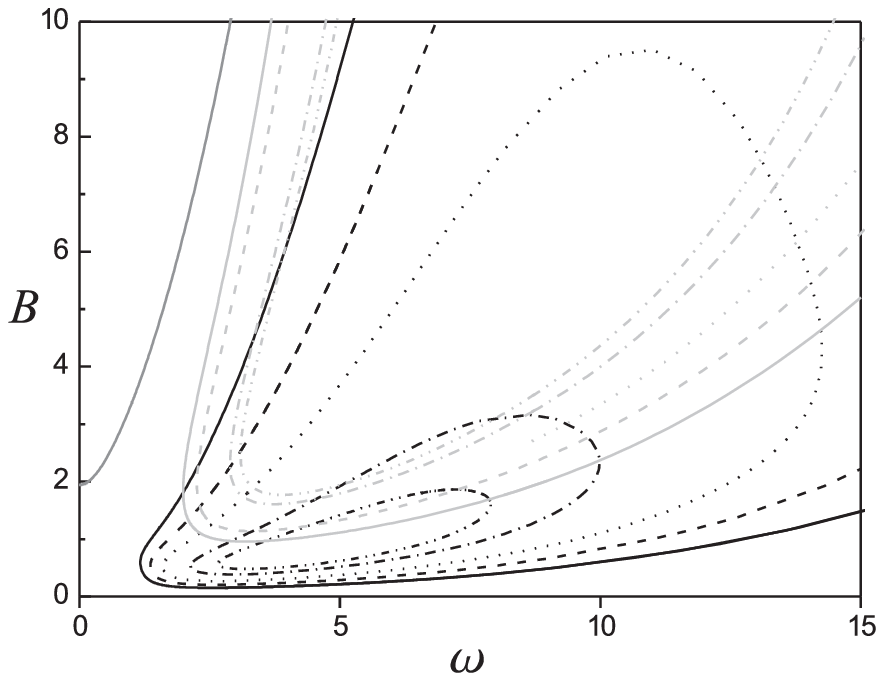}
\end{center}
\caption{Curves of secondary pitchfork bifurcations (dark grey),
Hopf bifurcations of $\lambda_1$ (black) and Hopf bifurcations of
$\lambda_2$ (light grey) at $\gamma=0.01$, 0.013, 0.016, 0.02 and
0.022 (from outside to inside). The stable region is below the dark
grey curve and around the black and light grey knees.}
\label{fig:hopf}
\end{figure}

Next we consider the combined effect of both damping and spin on the
stability of the rod. Figure~\ref{fig:tet:helix_stab-b1_B35_gamma05_w-real}
and \ref{fig:tet:helix_stab-b1_B35_gamma05_w-imag} show the evolution of the
first three eigenvalues as $\omega$ is increased for a solution along the
first bifurcating branch b1 at constant $B=1.2$ and two different values of
the damping coefficient. The value $\gamma = 0.04375$ was estimated by
Valverde et al.~\cite{valverde-cosserat} to be the damping value for a real
electrodynamic tether under working conditions. The figures show that
the solution remains stable at least up to $\omega=4$ (higher-order
eigenvalues do not seem to alter this picture). However, at the lower value
of $\gamma=0.01$ we observe a Hopf bifurcation at $\omega\approx 1.5$ where
the conjugate pair of eigenvalues ($\lambda_1$, $\bar\lambda_1$) acquires a
positive real part (at $\omega\approx 2.2$ the second pair of eigenvalues
($\lambda_2$, $\bar\lambda_2$) moves into the positive real half plane,
returning to the negative half plane at $\omega\approx 6$). Also note that
damping only affects the real parts of the eigenvalues, not the imaginary
parts (Figure~\ref{fig:tet:helix_stab-b1_B35_gamma05_w-imag}).

Figure \ref{fig:hopf} shows curves of these Hopf bifurcations, both
of $\lambda_1$ and $\lambda_2$, in the $B$-$\omega$ parameter plane
for various values of $\gamma$, as well as the curve of secondary
pitchfork bifurcations. The latter curve is independent of $\gamma$.
This makes sense as our (quasi-) statics problem, and hence the
bifurcation diagrams in Figures~\ref{fig:tet:bif-w025-measure1},
\ref{fig:tet:bif-w075-measure1} and \ref{fig:tet:bif-w2-measure1},
do not depend on $\gamma$, and therefore the secondary pitchfork
bifurcations, where a real eigenvalue goes through zero, must also
be independent of $\gamma$. Interestingly, the Hopf bifurcation
curves of $\lambda_1$ are closed, with instability inside the
curves, while the Hopf bifurcation curves of $\lambda_2$ are open.
The closed curves shrink to zero at $\gamma=0.023$, leaving only
instability due to Hopf bifurcations of $\lambda_2$ for larger
values of $\gamma$. In the limit $\gamma\to 0$ the Hopf curves
approach the coordinate axes $\omega=0$, $B=0$. The region of stable
b1 solutions in the diagram is the connected L-shaped region between
the secondary pitchfork and Hopf curves.

\section{Conclusion}

We have shown that whirling current-carrying rods bifurcate under
increasing magnetic field (or current). For the welded boundary conditions
considered here this bifurcation is doubly degenerate (even for a
transversely isotropic rod). This is because the equations are invariant
under rotation about the axis ($\bo e_3$ axis) of the supports. This
symmetry property complicated Wolfe's analysis, both for the whirling string
\cite{wolfe2} and the statical rod \cite{wolfe3}. In both models account had
to be taken of the variational nature of the problem to prove existence of
non-trivial bifurcating states. We find that for the rod, eigenvalues become
simple when rotating states are considered, provided the rod is anisotropic.
Neither spin nor anisotropy alone resolves the degeneracy.

Wolfe's analysis \cite{wolfe1} for non-whirling strings suggests
that the first bifurcating branch is stable and that all other
branches are unstable. We find the same to be true for a welded rod,
initially, but a secondary pitchfork bifurcation occurs where the
post-buckling solution (along b1) loses stability. Unfortunately, our method
allows us only to study stationary or quasi-stationary (whirling) solutions,
so it is not clear what type of solutions occur when the (quasi-) stationary
solutions lose stability in these secondary bifurcations. To investigate this
one would have to do simulations based on direct discretisation of the PDEs
(\ref{eq:lmomentum4}) and (\ref{eq:amomentum4}).

Associated with the secondary pitchfork bifurcations are connecting orbits
(b12, b34, etc.) which interfere with the primary pitchfork bifurcations
in the limit $\omega\to 0$. Another interesting bifurcation scenario concerns
the merging of symmetrically related branches (at $B=0$) and the subsequent
'lifting off' of these branches, thereby becoming non-bifurcating. In the
case of b1 (Figure~\ref{fig:tet:bif-w075-measure1}) this leads to the
disappearance of stable trivial (straight) solutions.

We have performed the, as far as we are aware, first study of the
magnetically induced post-buckling behaviour of elastic rods. Our results
may be of interest to the design of static as well as spinning electrodynamic
tethers, which often operate in the post-buckling regime. As we pointed out
in Section 7.1, the critical values of the magnetic field (or current) that
we find may be close to operating values for certain tether types.
Furthermore, axial spin tends to destabilise the post-buckling state of
the tether through Hopf bifurcations, although typical viscoelastic material
damping levels go some way towards avoiding Hopf instabilities.
On the other hand, spin has a stabilising effect in that it tends to push the
secondary pitchfork bifurcation to higher values of the magnetic field.
We have mapped out stable regions in the $B$-$\omega$ parameter plane
(Figure~\ref{fig:hopf}).


We finally like to speculate on another possible application of this
work. There is great current interest in conducting nanowires. These
can either be silicon-based wires, carbon nanotubes or metal-coated
biological fibres such as proteins, DNA molecules and microtubules
\cite{scheibel_et_al}. In addition there is the ongoing discussion
whether or not DNA molecules are electrical conductors
\cite{maiya_ramasarma}. All these nanometer-scale structures are
believed to have great potential as building blocks for future
electronic devices. The interaction of these wires with magnetic
fields could conceivably be exploited to obtain certain desirable
properties.



%


\newpage

\appendix
\begin{center}
\textbf{Appendix A: Matrices for the linearisation}
\end{center}

%

\noindent
The matrices $\bo B_i$ appearing in equation (\ref{per:lmomentum}) are given
by

\begin{equation}
\bo B_1 = \left(\begin{array}{ccc}
 0& -\kappa_3^0 & \kappa_2^0\\
 \kappa_3^0& 0 & -\kappa_1^0\\
 -\kappa_2^0& \kappa_1^0 &0
\end{array}\right),
\nonumber
\end{equation}

\begin{equation}
\bo B_2 = \omega^2 \left(\begin{array}{ccc}
 d_{11}^0& d_{12}^0 & 0\\
 d_{21}^0& d_{22}^0 & 0\\
 d_{31}^0& d_{32}^0 & 0
\end{array}\right),
\nonumber
\end{equation}

\begin{equation}
\bo B_3 =  \left(\begin{array}{ccc}
 0& F_{3}^0 & -F_{2}^0\\
 -F_{3}^0& 0 & F_{1}^0\\
 F_{2}^0& -F_{1}^0 & 0
\end{array}\right),
\nonumber
\end{equation}

\begin{equation}
\bo B_4 =  \left(\begin{array}{ccc}
 F_2^0\kappa_2^0+F_3^0\kappa_3^0-B(d_{22}^0d_{11}^0-d_{21}^0d_{12}^0)& (F_3^0)'-F_1^0\kappa_2^0 & -(F_2^0)'-F_1^0\kappa_3^0\\
 -(F_3^0)'-F_2^0\kappa_1^0& F_3\kappa_3^0+F_1^0\kappa_1^0-B(d_{22}^0d_{11}^0-d_{21}^0d_{12}^0) & (F_1^0)'-F_2^0\kappa_3^0\\
 (F_2^0)'-F_3\kappa_1^0-B(d_{22}^0d_{31}^0-d_{21}^0d_{32}^0)& -(F_1^0)'-F_3^0\kappa_2^0+B(d_{12}^0d_{31}^0-d_{11}^0d_{32}^0) &
 F_1^0\kappa_1^0+F_2^0\kappa_2^0
\end{array}\right),
\nonumber
\end{equation}

\begin{equation}
\bo B_5 = \left(\begin{array}{ccc}
 d_{11}^0& d_{12}^0 & d_{13}^0\\
 d_{21}^0& d_{22}^0 & d_{23}^0\\
 d_{31}^0& d_{32}^0 & d_{33}^0
\end{array}\right),
\nonumber
\end{equation}

\begin{equation}
\bo B_6 = 2\omega \left(\begin{array}{ccc}
 d_{12}^0& -d_{11}^0 & 0\\
 d_{22}^0& -d_{21}^0 & 0\\
 d_{32}^0& -d_{31}^0 & 0
\end{array}\right).
\nonumber
\end{equation}

Matrices $\bo C_i$ appearing in equation (\ref{per:amomentum}) are given by

\begin{equation}
\bo C_1 = \bo B_1, \nonumber
\end{equation}

\begin{equation}
\bo C_2 = \left(\begin{array}{ccc}
 0& M_{3}^0 & -M_{2}^0\\
 -M_{3}^0& 0 & M_{1}^0\\
 M_{2}^0& -M_{1}^0 & 0
\end{array}\right),
\nonumber
\end{equation}

\begin{equation}
\bo C_3 = \left(\begin{array}{ccc}
 C_3^{11}& C_3^{12} & C_3^{13}\\
 C_3^{21}& C_3^{22} & C_3^{23}\\
 C_3^{31}& C_3^{32} & C_3^{33}
\end{array}\right),
\nonumber
\end{equation}
where

\begin{eqnarray}
C_3^{11}&=&M_3^0\kappa_3^0+M_2^0\kappa_2^0-P(\bo \omega \cdot \bo d_3^0)(\bo d_2^0 \times \bo \omega \cdot \bo d_1^0)-P(\bo \omega \cdot \bo d_2^0)(\bo d_3^0 \times \bo \omega \cdot \bo d_1^0), \nonumber \\
C_3^{12}&=&(M_3^0)'-M_1^0\kappa_2^0-PR(\bo \omega \cdot \bo d_1^0)(\bo d_1^0 \times \bo \omega \cdot \bo d_3^0), \nonumber\\
C_3^{13}&=&-(M_2^0)'-M_1^0\kappa_3^0-F_1^0+PR(\bo \omega \cdot \bo d_1^0)(\bo d_1^0 \times \bo \omega \cdot \bo d_2^0)+P(\bo \omega \cdot \bo d_1^0)(\bo d_2^0 \times \bo \omega \cdot \bo d_1^0), \nonumber\\
C_3^{21}&=&-(M_3^0)'-M_2^0\kappa_1^0+P(\bo \omega \cdot \bo d_2^0)(\bo d_2^0 \times \bo \omega \cdot \bo d_3^0), \nonumber\\
C_3^{22}&=&M_3^0\kappa_3^0+M_1^0\kappa_1^0+PR(\bo \omega \cdot \bo d_3^0)(\bo d_1^0 \times \bo \omega \cdot \bo d_2^0)+PR(\bo \omega \cdot \bo d_1^0)(\bo d_3^0 \times \bo \omega \cdot \bo d_2^0), \nonumber \\
C_3^{23}&=&(M_1^0)'-M_2^0\kappa_3^0-F_2^0-PR(\bo \omega \cdot \bo d_2^0)(\bo d_1^0 \times \bo \omega \cdot \bo d_2^0)-P(\bo \omega \cdot \bo d_2^0)(\bo d_2^0 \times \bo \omega \cdot \bo d_1^0), \nonumber\\
C_3^{31}&=&(M_2^0)'-M_3^0\kappa_1^0+F_1^0-P(\bo \omega \cdot \bo d_3^0)(\bo d_2^0 \times \bo \omega \cdot \bo d_3^0), \nonumber\\
C_3^{32}&=&-(M_1^0)'-M_3^0\kappa_2^0+F_2^0+PR(\bo \omega \cdot \bo d_3^0)(\bo d_1^0 \times \bo \omega \cdot \bo d_3^0), \nonumber\\
C_3^{33}&=&M_2^0\kappa_2^0+M_1^0\kappa_1^0+P(1-R)(\bo \omega \cdot \bo d_2^0)(\bo d_1^0 \times \bo \omega \cdot \bo d_3^0)+P(1-R)(\bo \omega \cdot \bo d_1^0)(\bo d_2^0 \times \bo \omega \cdot \bo d_3^0),\nonumber \\
\nonumber
\end{eqnarray}

\begin{equation}
\bo C_4 = \left(\begin{array}{ccc}
 0& -1 & 0\\
 1& 0 & 0\\
 0& 0 & 0
\end{array}\right),
\nonumber
\end{equation}

\begin{equation}
\bo C_5 = P\left(\begin{array}{ccc}
 1& 0 & 0\\
 0& R & 0\\
 0& 0 & (1+R)
\end{array}\right),
\nonumber
\end{equation}

\begin{equation}
\bo C_6 = 2P \left(\begin{array}{ccc}
 (\bo d_2^0 \times \omega \times \bo d_3^0) \cdot \bo d_1^0& 0 & -(\bo d_2^0 \times \omega \times \bo d_1^0) \cdot \bo d_1^0\\
 0& -R(\bo d_1^0 \times \omega \times \bo d_3^0) \cdot \bo d_2^0 & R(\bo d_1^0 \times \omega \times \bo d_2^0) \cdot \bo d_2^0\\
 (\bo d_2^0 \times \omega \times \bo d_3^0) \cdot \bo d_3^0& -R(\bo d_1^0 \times \omega \times \bo d_3^0) \cdot \bo d_3^0 & R(\bo d_1^0 \times \omega \times \bo d_2^0) \cdot \bo
 d_3^0-(\bo d_2^0 \times \omega \times \bo d_1^0) \cdot \bo d_3^0
\end{array}\right).
\nonumber
\end{equation}

Matrices $\bo D_i$ appearing in equation (\ref{per:constitutive}) are given by

\begin{equation}
\bo D_1 = \frac{1}{f}\left(\begin{array}{ccc}
 -1& 0 & 0\\
 0& -R & 0\\
 0& 0 & -\frac{\Gamma(1+R)}{2}
\end{array}\right),
\nonumber
\end{equation}

\begin{equation}
\bo D_2 = \left(\begin{array}{ccc}
 0& \frac{1}{f}\kappa_3^0 & -M_2^0-\frac{1}{f}(1-R)\kappa_2^0\\
 -\frac{1}{f}\kappa_3^0& 0 & M_1^0-\frac{1}{f}(1-R)\kappa_1^0\\
 M_2^0-\frac{1}{f}\kappa_2^0\left(R-\frac{\Gamma(1+R)}{2}\right)&
 -M_1^0+\frac{1}{f}\kappa_1^0\left(1-\frac{\Gamma(1+R)}{2}\right) &0
\end{array}\right),
\nonumber
\end{equation}

\begin{equation}
\bo D_3 = \frac{\gamma}{f}\left(\begin{array}{ccc}
 0& -\kappa_3^0 & \kappa_2^0\\
 R\kappa_3^0& 0 & -R\kappa_1^0\\
 -\frac{\Gamma(1+R)}{2}\kappa_2^0& \frac{\Gamma(1+R)}{2}\kappa_1^0 &0
\end{array}\right),
\nonumber
\end{equation}

\begin{equation}
\bo D_4 = -\gamma \bo D_1. \nonumber
\end{equation}

All the $\kappa_i^0$ in the above can be expressed in terms of the moments
$M_i^0$ by means of the constitutive relations (\ref{stat:constitutive}).

\end{document}